\documentclass[english,reqno]{amsart}

\usepackage{diagbox}
\usepackage{amsmath, appendix, ulem, amsthm}
\usepackage[nobysame]{amsrefs}
\usepackage{amssymb, xcolor}
\usepackage[margin=1in]{geometry}
\usepackage{mathrsfs}
\usepackage{graphicx}
\usepackage{subfig}
\usepackage{float}
\usepackage{epsf}
\usepackage{titletoc}
\usepackage{cases}
\usepackage[linesnumbered,ruled]{algorithm2e}
\usepackage{colortbl}

\numberwithin{equation}{section}

\newtheorem{thm}{Main Theorem}
\theoremstyle{remark}

\newtheorem{remark}[thm]{Remark}

\renewcommand{\tilde}{\widetilde}

\newcommand{\f}{\frac}
\newcommand{\p}{\partial}
\newcommand{\beq}{\begin{equation}}
\newcommand{\eeq}{\end{equation}}

\newcommand{\dv}{ \, {\rm d} v}

\newcommand{\eps}{\varepsilon}

\newcommand{\RR}{\mathbb{R}}

\newcommand{\rd}{\mathrm{d}}

\newcommand{\vpran}[1]{\left(#1\right)}

\newcommand{\Dt}{\Delta t}
\newcommand{\Dx}{\Delta x}
\newcommand{\Uns}{U^{n+1*}}
\newcommand{\Ins}{I^{n+1*}}
\newcommand{\half}{\frac{1}{2}}
\newcommand{\order}{\mathcal{O}}
\newcommand{\average}[1]{\left\langle#1\right\rangle}
\newcommand{\dc}{D_d}
\newcommand{\ihalf}{{i+\frac{1}{2}}}
\newcommand{\nps}{{n+1*}}

\newcommand{\Vis}{K_i^{n+1*}}
\newcommand{\Vips}{K_{i+\half}^{n+1*}}

\newcommand{\black}[1]{\textcolor{black}{#1}}

\begin{document}

\title{Accurate front capturing asymptotic preserving scheme for nonlinear gray radiative transfer equation}
\author{Min Tang}
\address{School of Mathematics, Institute of Natural Sciences and MOE-LSC, Shanghai Jiao Tong University, China}
\email{tangmin@sjtu.edu.cn}
\author{Li Wang}
\address{School of Mathematics, University of Minnesota, 206 Church St SE, Minneapolis, MN 55455}
\email{wang8818@umn.edu}
\author{Xiaojiang Zhang}
\address{School of Mathematics, Institute of Natural Sciences, Shanghai Jiao Tong University, China}
\email{xjzhang123@sjtu.edu.cn}
\thanks{M.Tang is supported by Science Challenge Project No. TZ2016002 and NSFC11871340 and NSFC91330203. L.Wang is partially supported by NSF-DMS 1903420, NSF CAREER-DMS 1846854 and a start-up fund from University of Minnesota.}
\date{\today}
\begin{abstract}
\textcolor{black}{We develop an asymptotic preserving scheme for the gray radiative transfer equation. Two asymptotic regimes are considered: one is a diffusive regime described by a nonlinear diffusion equation for the material temperature; the other is a free streaming regime with zero opacity.} To alleviate the restriction on time step and capture the correct front propagation in the diffusion limit, an implicit treatment is crucial. However, this often involves a large-scale nonlinear iterative solver as the spatial and angular dimensions are coupled. Our idea is to introduce an auxiliary variable that leads to a ``redundant" system, which is then solved with a three-stage update: prediction, correction, and projection. The benefit of this approach is that the implicit system is local to each spatial element, independent of angular variable, and thus only requires a scalar Newton's solver. We also introduce a spatial discretization with a compact stencil based on even-odd decomposition. \textcolor{black}{Our method preserves both the nonlinear diffusion limit with correct front propagation speed and the free streaming limit, with a hyperbolic CFL condition.}
\end{abstract}

\maketitle

\section{Introduction}

The gray radiative transfer equation (GRTE) concerns photon transport and its interaction with the 
background material. It describes the radiative transfer and energy exchange between radiation and materials, and has wide applications in astrophysics and inertial confinement fusion.
The system for the radiative intensity $I$ and the material temperature $T$ is
\begin{equation}
\left\{
\begin{aligned}
& \frac{1}{c} \partial_t I+\frac{1}{\varepsilon} v\cdot \nabla_x I=
\frac{\sigma}{\varepsilon^2}\left(acT^4-I\right), \quad x  \in \Omega_x \subset \RR^d, ~ v\in \Omega_v \subset \RR^d \,;\\
&  C_v\partial_t T=\frac{\sigma}{\varepsilon^2}(\rho-acT^4), \qquad \rho=\f{1}{|V|}\int_{V}I \dv := \average{\rho}\,,
\end{aligned}
\right.
\label{radiative_transport}
\end{equation}
which is equipped with initial and boundary condition 
\begin{equation*}
I(0,x,v) = I_{\text{in}}(x,v), \quad T(0,x) = T_{\text{in}}(x)\,,
\end{equation*}
and
\begin{equation} \label{BC00}
I(t, x, v) = b(t,v), \quad v\cdot n_x <0\,.
\end{equation}
In equation \eqref{radiative_transport}, $x$ is the space variable, $v$ is the angular variable, $\sigma(x,T)$ is the opacity and $a$, $c$, $C_v$ are positive constants
representing the radiation constant, light speed and heat capacity respectively.
\textcolor{black}{The parameter $\varepsilon>0$ is the scaled mean free path that describes the distance of two successive scatterings} and $n_x$ in \eqref{BC00} is the outer normal direction for $x \in \partial \Omega_x$.

Due to the high dimensionality and nonlinearity in \eqref{radiative_transport}, it is often expensive to directly discretize \eqref{radiative_transport}. Therefore, some limit cases are identified. \textcolor{black}{Two limit cases are of particular interest, one is the free streaming case, in which $\sigma=0$ and the photons transport freely without change of velocity. The other is when the mean free path is small and the unknowns vary slowly in time and space, i.e. the parabolic scaling limit by sending $\varepsilon\to 0$.} In the parabolic scaling, the radiative intensity $I$ approaches
a Planckian at the local temperature 
$acT^4$ and the  
material temperature satisfies the following nonlinear diffusion equation \cite{LPB,BGP}: 
\begin{equation}\label{eq:limit_T}
a \partial_t T^4+C_v \partial_t T=\nabla_x \cdot \vpran{\frac{ac}{\sigma \dc }\nabla_x T^4}\,,
\end{equation}
where $\dc$ is a constant depending on the dimension of the velocity space. For instance $\dc = 3$ in one dimensional case. To derive the limiting equation in \eqref{eq:limit_T}, one can add up the two equations in \eqref{radiative_transport} and consider the system
\begin{equation*}
\left\{
\begin{aligned}
& \frac{1}{c} \partial_t I+C_v\partial_t T+\frac{1}{\eps} v\cdot \nabla_x I=
\frac{\sigma}{\eps^2}\left(\rho-I\right)\,; \\
&  C_v\partial_t T=\frac{\sigma}{\eps^2}(\rho-acT^4).
\end{aligned}
\right.
\end{equation*}
Then by using Chapman-Enskog expansion, the leading order temperature satisfies the nonlinear diffusion equation \eqref{eq:limit_T}.

The design of an efficient 
scheme for \eqref{radiative_transport} has a three-fold difficulty. On the one hand, for high opacity material when the interactions between radiation and material are strong, the mean free path is very small, which indicates that $\varepsilon\ll 1$. In order to solve the GRTE accurately, one has to use resolved space and time steps that are less than the mean free path, which leads to an extremely high computational cost. One way around it is to simulate the nonlinear diffusion equation \eqref{eq:limit_T} instead. However, when the opacity $\sigma$ depends on $T$, both optical thin and optical thick regimes co-exist, solving only the limit model will not generate satisfactory results. 
\black{On the other hand, the nonlinearity induces additional complexity. The time step of explicit schemes is limited by a CFL condition based on the speed of light, thus the Knudsen number $\varepsilon$. In order to allow larger time steps, implicit schemes have to be employed, then a nonlinear system has to be solved for each time step. Assume that there are respectively $N_x$ and $N_v$ grids in space and velocity. Since $I$ and $T$ are coupled together,  a $(N_xN_v+N_x)\times (N_xN_v+N_x)$ linear system has to be solved for each iteration of the nonlinear solver, which is expensive. 
Furthermore, the degeneracy in the limit equation calls for special treatment. Indeed, the nonlinear diffusion limit \eqref{eq:limit_T} is degenerate at $T=0$, thus starting from a compactly supported initial intensity, the solution of the limit equation remains compactly supported with finite propagation speed. See the detailed discussion in \cite{LTWZ18} and reference therein. It is desirable to capture this speed numerically. However, as pointed out in \cite{Lowrie}, a nonlinear implicit solver is expected to serve the purpose, only if the linearization adopted in the nonlinear iterative solver is properly chosen. That is, a stable and convergent linearization may not lead to a correct front capturing method. }
 
Due to the above mentioned difficulties, designing efficient numerical solvers for GRTE is a major endeavor and remains a challenging task. To resolve the difficulty stemming from small scales, the asymptotic preserving (AP) schemes that provide seamless connections between thin and thick opacity regimes would serve the purpose. To name a few, such works include \cite{LM89, JL93, GJL99, Klar98,JTH09,HTY14} for the linear steady state RTE, \cite{JPT00, LM08, Qiu14, KFJ16} for time dependent problems, and \cite{BPR13} for high order methods. However, as explicit time dependent solvers are often limited either by the parabolic CFL condition or hyperbolic CFL condition based on the speed of light, implicit methods become the new paradigm. As a result, semi-implicit or fully implicit solvers that involve linearization or the development of pre-conditioner start to cumulate. For linear RTE, three kinds of methods prevail. One is the Krylov iterative method for the discrete-ordinate system preconditioned by diffusion synthetic acceleration (DSA) \cite{WWM04, LM83}, another is the implicit Monte Carlo (IMC) method \cite{FC71}, and the third is moment method in terms of spherical harmonics \cite{Pomraning93, Larsen80, BL01, FS11, Hauck11, OLS13}. Practically however, they all come with caveat. The first kind includes a sub-iterations of sweeps and diffusion solvers for DSA pre-conditioner, and thus is very complicated to extend to nonlinear case. The second kind suffers from unavoidable statistical fluctuations and requires a large sampling of points. The third one often generate very complicated system whose well-posedness and realizability necessitate detailed investigation. A more recently developed method that takes advantage of the spectral structure of the linear system resulting from a fully implicit discretization provides an viable alternative \cite{LW17}. For the nonlinear RTE, much less works are available. \textcolor{black}{  McClarren et. al. made a first attempt in designing a semi-implicit time integration that treats the streaming term explicitly and material coupling term implicitly with a linearized emission source, along with a moment method for angular variable and a discontinuous Galerkin method for spatial variables \cite{MELD08}. Another method is the unified gas kinetic scheme developed by Sun et.al., which employs a linearized iterative solver for the material thermal temperature \cite{SJX1,SJX2}. These  discretizations are AP in space, but their time steps are limited by a CFL condition based on the speed of light, thus the Knudsen number $\varepsilon$, while this restriction relaxes to explicit diffusion CFL condition in the diffusion limit \cite{Klar98,MELD08,SJX1,SJX2}. In order to allow larger time steps, an almost fully implicit time discretization was proposed in \cite{LKM13}. It treats all terms except the opacity implicit. The authors proved that the scheme satisfies maximum principle, thus can provide good result with large time step. However it calls for nonlinear iterative solver that can be quite expensive in high dimensions. All the above mentioned schemes either require $\Delta t\sim \max\{O(\varepsilon \Delta x),O(\Delta x^2)\}$ or nonlinear iteration for a large system.} In this paper, we aim at developing an AP scheme that requires only a {\it scalar} Newton's iterative solver \black{and allows a larger time step than the semi-implicit method in \cite{MELD08}.}

When the opacity nonlinearly depends on $T$, scattering coefficient introduces another level of nonlinearity. One particular example is $\sigma_T=\frac{\sigma}{T^3}$ with $\sigma$ independent of $T$. \black{This case is of particular interest, especially in the Marshak  wave simulations. Marshak wave is the nonlinear energy wave created when a cold material is heated at one end. Since the initial temperature is low and the final temperature is high, the problem is initially diffusive and becomes non-diffusive in the region through which the wavefront has passed.} 
 Rewriting \eqref{radiative_transport} to incorporate the temperature dependence in $\sigma$, 
\begin{equation}
\left\{
\begin{aligned}
& \frac{1}{c} \partial_t I+\frac{1}{\varepsilon} v\cdot \nabla_x I=
\frac{1}{\varepsilon^2} \frac{\sigma}{T^3}\left(acT^4-I\right), \quad x  \in \Omega_x \subset \RR^d, ~ v\in \Omega_v \subset \RR^d \,; \\
&  C_v\partial_t T=\frac{1}{\varepsilon^2} \frac{\sigma}{T^3}(\rho-acT^4),
\end{aligned}
\right.
\label{marshak00}
\end{equation}
then sending $\eps \rightarrow 0$, its diffusion limit reads
\begin{equation}\label{limit_marshak}
a \partial_t T^4+C_v \partial_t T=\nabla_x \cdot \vpran{\frac{ac T^3}{\sigma \dc }\nabla_x T^4}\,.
\end{equation}
The additional nonlinearity in the \black{opacity} calls for special care in the design of of nonlinear solver in order to capture the correct propagation speed of the front.

We propose a new scheme that can solve the above-mentioned difficulties, for both temperature dependent and independent opacity. The benefits of our scheme are: 
\begin{itemize}
\item The scheme is AP \black{in both diffusive and free streaming regimes.} Thus it allows to use unresolved meshes; and its stability and accuracy is independent of $\varepsilon$ whose magnitude ranges from very small to order one. 
\item \black{It can use a hyperbolic CFL condition $\Delta t = O(\Delta x)$ that is independent of $\varepsilon$ to
provide the correct solution behavior. On the contrary, the CFL numbers of previous AP method developed in \cite{MELD08,SJX1} depend on the unscaled light speed, thus the Knudsen number $\varepsilon$. } 
\item \black{The prediction step only requires a linear solver. Updating the material temperature $T$ and $I$ in the correction step are decoupled. $T$ requires only {\it scalar} Newton iteration for solving a forth order polynomial on each space grid, and 
 $I$ that depends both on space and angular variables requires only a {\it linear} solver.}
 
 \item It can capture the correct front position when the initial temperature is compact supported, \black{even when the opacity nonlinearly depends on $T$.} 
\item It can be extended to nonlinearities other than $T^4$, thus is very attractive for extensions to multi-frequency radiative transfer problems.   
\end{itemize}

The organization of this paper is as follows. Section 2 exposes a three-stage diffusion solver for the limiting nonlinear diffusion equation, it requires only {\it scalar} iterative solver and provides correct front speeds for compact supported initial data. 
The detailed semi-implicit
scheme and
one dimensional fully discretized schemes for GRTE are presented in Sections 3 and 4 respectively. Their capability to capture the nonlinear diffusion limit is shown as well.
In Section 5, how to deal with the case when $\sigma$ depends on $T$ is discussed and finally, some numerical tests are given in Section 6 to discuss
the stability and accuracy of our scheme. Both cases when $\sigma$ is with or without $T$ dependence are considered, with particular attentions been paid to the compactly supported initial conditions. Finally the paper is concluded in Section 7. 

\section{Three-stage diffusion solver}
Before treating the GRTE, we first discuss the discretization for the limiting nonlinear degenerate diffusion equation \eqref{eq:limit_T}. We intend to propose a scheme that (1) requires only hyperbolic CFL condition, (2) can avoid to solve a nonlinear {\it system}, (3) can capture the correct propagation front.
Our idea is to introduce an auxiliary variable whose equation along with \eqref{eq:limit_T} forms a ``redundant" system which mitigates the nonlinearity and can be solved semi-implicitly. This idea shares the same spirits as that in our previous work \cite{LTWZ18}, in which we proposed an accurate front capturing scheme for tumor growth models with similar nonlinearity. The details are described below.

\subsection{The prediction-correction-projection method}
Let $U=T^4$, then rewrite \eqref{eq:limit_T} into
$$
4a T^3\p_t T+C_v\p_t T=\nabla_x \cdot \vpran{\frac{ac}{\sigma \dc}\nabla_x U},
$$
which yields
\beq\label{eq:TT3}
\p_t T=\frac{1}{4aT^3+C_v} \nabla_x \cdot  \vpran{\frac{ac}{\sigma \dc}\nabla_x U} .
\eeq
After multiplying both sides of \eqref{eq:TT3} by $4T^3$, one gets a {\it semi-linear} equation for the new variable $U$ as follows
\beq\label{eq:limitU}
\p_t U=\frac{4T^3}{4aT^3+C_v} \nabla_x \cdot \vpran{\frac{ac}{\sigma \dc}\nabla_x U}.
\eeq
When the solution is regular enough, equation \eqref{eq:limitU} for $U$ is equivalent to equation \eqref{eq:limit_T} for $T$. 
Similar to the idea in \cite{LTWZ18}, we solve the two equations for $T$ and $U$ together,
\beq\label{eq:limit_TV}
\left\{\begin{aligned}
&\p_t U=\frac{4T^3}{4aT^3+C_v} \nabla_x  \cdot \vpran{\frac{ac}{\sigma \dc}\nabla_x U} ;\\
&a\partial_t T^4+C_v \partial_t T=\nabla_x \cdot \vpran{\frac{ac}{\sigma \dc}\nabla_x U},\\
\end{aligned}
\right.
\eeq
with the initial condition for $U$ as $U(x,0)=T^4(x,0)$.
Note specifically that, when solving \eqref{eq:limit_TV} numerically, the relation $U= T^4$ that bridges these two equations is polluted by the numerical error and therefore a projection step is needed to reinforce this relation. More specifically, we have the following semi-discretized prediction-correction-projection method.

At every time step, given $T^n$, we first update $U$ by a prediction step:
\beq\label{eq:limit_U_dis1}
\frac{U^{n+1\ast}-U^n}{\Delta t}=\frac{4(T^n)^3}{4a(T^{n})^3+C_v} \nabla_x \cdot \vpran{\frac{ac}{\sigma \dc}\nabla_x U^{n+1\ast}} .
\eeq
and then using the updated $U^{n+1\ast}$ to advance $T$ via \eqref{eq:limit_TV} as:
\begin{equation}\label{eq:limit_UT_dis2}
a\frac{(T^{n+1})^4-(T^n)^4}{\Delta t}+C_v \frac{T^{n+1}-T^n}{\Delta t}=\nabla_x \cdot \vpran{\frac{ac}{\sigma \dc}\nabla_x U^{n+1 \ast}}.
\end{equation}
Finally a projection is conducted to get $U^{n+1}$ by
\beq\label{eq:limit_UT_dis3}
U^{n+1}=(T^{n+1})^4.
\eeq

\subsection{Spatial discretization}
For simplicity, we consider 1D case. Let $\Delta x=L/N_x$, we consider uniform mesh as follows
$$x_i=(i-1)\Delta x,\qquad \mbox{for $i=1,\cdots, N_x+1,$}$$
and let 
\[
x_{i+1/2}=(x_i+x_{i+1})/2, \quad i = 1, 2 \cdots, N_x. 
\]

Let $U_{i+\half}^n$, $T_{i+\half}^n$ to be the approximation of $U( t_n, x_\ihalf)$ and $T(t_n,x_\ihalf)$, then the fully discrete three-stage scheme for the diffusion limit reads
\begin{equation} \label{dis-diff}
\left\{\begin{aligned}
&\frac{U_\ihalf^{n+1*} - U_\ihalf^n}{\Delta t}  = \frac{4 (T_\ihalf^n)^3}{4a (T_\ihalf^n)^3 + C_v} \frac{1}{\Dx} \left[   \left(\frac{ac}{3\sigma}\right)_{i+1} \frac{ U_{i+\frac{3}{2}}^{n+1*} - U_\ihalf^{n+1*}}{\Dx } - \left(\frac{ac}{3\sigma}\right)_{i} \frac{ U_\ihalf^{n+1*} - U_{i-\half}^{n+1*}}{\Dx } \right] \,;
\\ & a \frac{(T^{n+1} _\ihalf)^4 - (T_\ihalf^n)^4}{\Dt} + C_v \frac{T_\ihalf^{n+1}- T_\ihalf^n}{\Dt }  =\frac{1}{\Dx} \left[   \left(\frac{ac}{3\sigma}\right)_{i+1} \frac{ U_{i + \frac{3}{2}}^{n+1*} - U_{i+\half}^{n+1*}}{\Dx } - \left(\frac{ac}{3\sigma}\right)_{i} \frac{ U_{i+\half}^{n+1*} - U_{i-\half}^{n+1*}}{\Dx } \right] \,;\\
& U_{\ihalf}^{n+1}=(T_{\ihalf}^{n+1})^4\,,
\end{aligned} \right. 
\end{equation}
for $1\leq i \leq N_x$. In order to be consistent with the discretization of GRTE in section 3, here we let the unknowns be on half grids.

\begin{remark}Due to the time derivative terms in the nonlinear diffusion limit for $T$, to construct a conservative scheme, it is unavoidable to solve nonlinear equations.
However, it is important to notice that, all nonlinear equations for different grid points are {\it decoupled} in the second stage of \eqref{dis-diff}. Thus a nonlinear equation of only one variable has to be solved for each grid point. One can use standard iterative solvers such as Newton's method for such polynomial equations and the convergence is guaranteed. 
\end{remark}

\begin{remark} It is necessary to use the three-stage scheme in order to meet all three requirements at the beginning of this section. To allow hyperbolic CFL condition, the diffusion operator can not be explicit, while if it is implicit, a nonlinear system that couples the unknowns at all grid points has to be solved. The introduction of the auxiliary variable $U$ provides a good approximation for $T^4$ by solving a linear system. However, the correct front speed can not be obtained by discretizing only the equation for $U$. Indeed, if we have only two stages---the prediction and projection: 
\begin{equation} \label{dis-diff-2}
\left\{\begin{aligned}
&\frac{U_\ihalf^{n+1} - U_\ihalf^n}{\Delta t}  = \frac{4 (T_\ihalf^n)^3}{4a (T_\ihalf^n)^3 + C_v} \frac{1}{\Dx} \left[   \left(\frac{ac}{3\sigma}\right)_{i+1} \frac{ U_{i+\frac{3}{2}}^{n+1} - U_\ihalf^{n+1}}{\Dx } - \left(\frac{ac}{3\sigma}\right)_{i} \frac{ U_\ihalf^{n+1} - U_{i-\half}^{n+1}}{\Dx } \right]\,;\\
& T_{\ihalf}^{n+1}=(U_{\ihalf}^{n+1})^{1/4},
\end{aligned} \right.
\end{equation}
then for $T^n_{\ihalf}=0$, $T_\ihalf^{n+1}$ remains zero, which indicates that the support of compactly supported initial data remains the same. \end{remark}

\subsection{\black{The nonlinearity in the opacity}}
\black{When the scattering coefficient nonlinearly depends on T, the above scheme can be applied with little variations. Hereafter, we consider a specific example with $\sigma_T =\frac{\sigma}{T^3}$,} then our prediction-correction-projection method becomes
\begin{equation} \label{marshak-diff2}
\left\{\begin{aligned}
&\frac{U^{n+1\ast}-U^n}{\Delta t}=\frac{4(T^n)^3}{4a(T^{n})^3+C_v} \nabla_x \cdot \vpran{\frac{ac(T^n)^3}{\sigma \dc}\nabla_x U^{n+1\ast}} \,; \\
&T^{n+1\ast}=(U^{n+1\ast})^{1/4},
\\ & a\frac{(T^{n+1})^4-(T^n)^4}{\Delta t}+C_v \frac{T^{n+1}-T^n}{\Delta t}=\nabla_x \cdot \vpran{\frac{ac(T^{n+1\ast})^3}{\sigma  \dc}\nabla_x U^{n+1 \ast}} \,; \\
& U^{n+1}=(T^{n+1})^4.
\end{aligned} \right.
\end{equation}
The corresponding fully discretized scheme writes
\begin{equation*}
\left\{\begin{aligned}
&\frac{U_\ihalf^{n+1*} - U_\ihalf^n}{\Delta t}  = \frac{4 (T_\ihalf^n)^3}{4a (T_\ihalf^n)^3 + C_v} \frac{1}{\Dx} \left[   \left(\frac{ac (T^n)^3}{3\sigma}\right)_{i+1} \frac{ U_{i+\frac{3}{2}}^{n+1*} - U_\ihalf^{n+1*}}{\Dx } - \left(\frac{ac (T^n)^3}{3\sigma}\right)_{i} \frac{ U_\ihalf^{n+1*} - U_{i-\half}^{n+1*}}{\Dx } \right] \,; \\
&T_{\ihalf}^{n+1*}=(U_{\ihalf}^{n+1*})^{1/4}\,;
\\ & a \frac{(T^{n+1} _\ihalf)^4 - (T_\ihalf^n)^4}{\Dt} + C_v \frac{T_\ihalf^{n+1}- T_\ihalf^n}{\Dt }  =\frac{1}{\Dx} \left[   \left(\frac{ac (T^{n+1*})^3}{3\sigma}\right)_{i+1} \frac{ U_{i + \frac{3}{2}}^{n+1*} - U_{i+\half}^{n+1*}}{\Dx } \right.
\\&\left.\hspace{7.5cm} - \left(\frac{ac (T^{n+1*})}{3\sigma}\right)_{i} \frac{ U_{i+\half}^{n+1*} - U_{i-\half}^{n+1*}}{\Dx } \right] \,; \\
& U_{\ihalf}^{n+1}=(T_{\ihalf}^{n+1})^4\,.
\end{aligned} \right. 
\end{equation*}

\section{Time discretization for the radiative transfer equation}
We start off by considering the time discretization of GRTE \eqref{radiative_transport}. In order to preserve the asymptotic limit \eqref{eq:limit_T} while keeping the computational complexity under control, we propose a scheme whose limit mimics the three-stage diffusion solver \eqref{eq:limit_U_dis1}--\eqref{eq:limit_UT_dis3}. 

\subsection{The prediction step}
As with \eqref{eq:limit_U_dis1}, we will solve for variable $U$ instead of $T$ in the prediction step. Multiplying the second equation in \eqref{radiative_transport} by $4T^3$ yields
\begin{subequations} \label{pre00}
\begin{numcases}{}
\frac{1}{c} \partial_t I+\frac{1}{\varepsilon} v\cdot \nabla_x I=
\frac{\sigma}{\varepsilon^2}\left(acU-I\right) \,; \label{pre00a}\\
C_v\partial_t U=4T^3\frac{\sigma}{\varepsilon^2}(\rho-acU).\label{pre00b}
\end{numcases}
\end{subequations}
Then at each time step, given $I^n$, $T^n$, $U^n$, we use the following semi-discrete scheme in time \begin{subequations} \label{semi-pre}
\begin{numcases}{}
\frac{1}{c} \frac{I^{n+1\ast}-I^n}{\Delta t}+\frac{1}{\varepsilon} v\cdot \nabla_x I^{n+1\ast}=
\frac{\sigma}{\varepsilon^2}\left(acU^{n+1\ast}-I^{n+1\ast}\right)\,;  \label{semi-pre1}\\
 C_v\frac{ U^{n+1\ast}-U^n}{\Delta t}=4\frac{\sigma}{\varepsilon^2}(T^n)^3(\rho^{n+1\ast} -acU^{n+1\ast}). \label{semi-pre2}
\end{numcases}
\label{radiative_IU}
\end{subequations}
Note that, considering $T$ given at previous time step, \eqref{pre00} is a linear system with respect to $I$ and $U$, and any asymptotic preserving scheme for linear transport equations such as micro-macro decomposition \cite{LM08} or even-odd decomposition can be used \cite{JPT00}. Here we use a fully implicit scheme to maximize the stability.
 
\subsection*{Diffusion limit of \eqref{radiative_IU}:} First from \eqref{semi-pre1}, one has
\begin{align}
I^{n+1*} &= ac \Uns - \frac{\eps^2}{\sigma } \left[ \frac{1}{c} \frac{I^{n+1*} - I^n}{\Dt }  + \frac{1}{\eps} v \cdot \nabla_x I^{n+1*}\right] \nonumber
\\ & =  ac \Uns - \frac{\eps}{\sigma} v \cdot \nabla_x (ac \Uns) + \order(\eps^2)  \,. \label{916}
\end{align}
Taking the average in $v$ for \eqref{semi-pre1}, dividing \eqref{semi-pre2} by $4 (T^n)^3$, and then adding these two equations, we have, to the leading order:
\begin{equation} \label{9161}
\frac{1}{c} \frac{\rho^{n+1*} - \rho^n}{\Dt}  + \average{v \cdot \nabla_x \left(  - \frac{1}{\sigma} v \cdot \nabla_x (ac \Uns) \right) } + \frac{C_v}{ 4 (T^n)^3} \frac{\Uns - U^n}{\Dt}  = 0\,,
\end{equation}
where we have plugged the expression of $I^{n+1*}$ from \eqref{916} into \eqref{semi-pre1}. Note from \eqref{semi-pre2} that 
\[
\rho^{n+1*} = ac \Uns + \order(\eps^2)\,,
\]
\eqref{9161} reduces to 
\begin{equation} \label{limit_AP_1}
\left[ a + \frac{C_v}{ 4 (T^n)^3} \right] \frac{\Uns - U^n}{\Dt} = \nabla_x \cdot \left(  \frac{ac}{\sigma \dc} \nabla_x \Uns\right)\,,
\end{equation}
which is the same as \eqref{eq:limit_U_dis1}. It is important to point out that \eqref{limit_AP_1} is only valid when $\rho^n = ac U^n+O(\varepsilon)$, and therefore a projection step as will be explained in Section \ref{sec:projection} is necessary.

\subsection{The correction step}
In order to mimic the correction stage \eqref{eq:limit_UT_dis2} when sending $\eps$ to zero, we first introduce the following expansion for $I$:
\black {\beq\label{J00}
I=ac T^4 +\varepsilon J\,,
\eeq
and substitute it into \eqref{radiative_transport}, which yields
\begin{equation}\label{radiative_transport3}
\left\{\begin{aligned}
& \f{1}{c}\p_t\left(ac T^4+ \varepsilon J\right)+\f{1}{\varepsilon}v\cdot \nabla_x I=-\f{\sigma}{\varepsilon} J \,;\\
& C_v\p_tT= \f{\sigma}{\varepsilon}\rho_J\,,
\end{aligned}
\right.
\end{equation}
where $\rho_J=\f{1}{|V|}\int_V J^{n+1}\dv$. Then, we try to solve \eqref{radiative_transport3} using the predicted information $\Uns$, $\Ins$. 
We use the following semi-discrete scheme in time
\begin{subequations} \label{semi-cor0}
\begin{numcases}{}
 \f{1}{c}\left(ac\f{(T^{n+1})^4-(T^{n})^4}{\Delta t}+\varepsilon \f{\tilde J^{n+1}-J^n}{\Delta t}\right)+\f{1}{\varepsilon}v\cdot \nabla_x  I^{n+1\ast}=-\f{\sigma}{\varepsilon} \tilde J^{n+1} \,; \label{semi-cor10}\\
 C_v\frac{T^{n+1}-T^n}{\Delta t}= \f{\sigma}{\varepsilon}\tilde \rho_J^{n+1}. \label{semi-cor20}
\end{numcases}
\end{subequations}
The advantage of treating the term $\f{1}{\varepsilon}v\cdot\nabla I$ using the information obtained from the prediction stage can be seen as follows. Taking the average of \eqref{semi-cor10} with respect to $v$, summing it up with \eqref{semi-cor20}, we have the following equations for $T^{n+1}$ and $\tilde\rho_J^{n+1}$:
\begin{subequations} \label{semi-TI}
\begin{numcases}{}
 \f{1}{c}\left(ac\f{(T^{n+1})^4-(T^n)^4}{\Delta t}+\varepsilon\f{\tilde\rho_J^{n+1}-\rho_J^n}{\Delta t}\right)+ C_v\f{T^{n+1}-T^n}{\Delta t}+ \f{1}{\varepsilon} \average{ v\cdot\nabla_x I^{n+1\ast}} =0; \label{semi-TI10}\\
 C_v\frac{T^{n+1}-T^n}{\Delta t}= \f{\sigma}{\varepsilon}\tilde\rho_J^{n+1}. \label{semi-TI20}
\end{numcases}
\end{subequations}
In particular, to update $T$, one only needs to solve $scalar$ nonlinear polynomial equations at each step, which is much more efficient than solving a nonlinear system. After getting $T^{n+1}$, $J^{n+1}$ can be updated implicitly by 
\beq\label{semi-J}
\f{1}{c}\left(ac\f{(T^{n+1})^4-(T^{n})^4}{\Delta t}+\varepsilon \f{J^{n+1}-J^n}{\Delta t}\right)+\f{1}{\varepsilon}v\cdot \nabla_x (ac(T^{n+1})^4+\varepsilon J^{n+1})=-\f{\sigma}{\varepsilon} J^{n+1}. 
\eeq
 }
\subsection*{Diffusion limit of \eqref{semi-cor0}:}
 \black{Note that, $$I^{n+1*} = ac \Uns - \frac{\eps}{\sigma} v \cdot \nabla_x (ac \Uns)+O(\eps^2)$$ from the prediction stage. Substitute it into \eqref{semi-TI10}, we have:
 \begin{equation*}
 \f{1}{c}\left(ac\f{(T^{n+1})^4-(T^n)^4}{\Delta t}+\varepsilon\f{\rho_J^{n+1}-\rho_J^n}{\Delta t}\right)+ C_v\f{T^{n+1}-T^n}{\Delta t}- \average{ v\cdot\nabla_x \left(\frac{1}{\sigma} v \cdot \nabla_x (ac \Uns)\right) }=O(\eps)\ .	
 \end{equation*}
Sending $\varepsilon\to 0$ in the above equation, formally one gets:
$$
a\f{(T^{n+1})^4-(T^n)^4}{\Delta t}+C_v\f{T^{n+1}-T^n}{\Delta t}=\nabla_x \cdot \left(\f{ac}{\sigma\dc}\nabla_x U^{n+1\ast}\right)\,,
$$
which is the same as \eqref{eq:limit_UT_dis2}.}

\subsection{The Projection step} \label{sec:projection}
Finally, we apply the projection step. \black{Given $T^{n+1}$ and $J^{n+1}$ from the correction step, we update $U^{n+1}$ and finally $I^{n+1}$ as follows:
\begin{equation} \label{proj00}
U^{n+1}=(T^{n+1})^4,\qquad  I^{n+1}=acU^{n+1}+\varepsilon J^{n+1}\,.
\end{equation}}
One sees that when $\varepsilon\to 0$, 
\[
I^{n+1}=acU^{n+1}=ac(T^{n+1})^4\,,
\]
which also indicates that to the leading order, $\average{I^{n+1}} = ac U^{n+1}$ for $\forall n >1$, and is consistent with \eqref{eq:limit_UT_dis3}. 

In summary, we have the following semi-discrete update at each time step:
 \begin{algorithm}[h]
\caption{one step of semi-discrete update for GRTE}
\SetAlgoLined
\KwIn{$I^n$, $U^n$, $T^n$, ($U^n = (T^{n})^4$) }
\KwOut{$I^{n+1}$, $U^{n+1}$, $T^{n+1}$, ($U^{n+1} = (T^{n+1})^4$) }

\BlankLine
    \textrm{\bf Prediction}: obtain $I^{n+1*}$, $U^{n+1*}$ from \eqref{semi-pre}; 
    \\ \textrm{\bf Correction}: obtain $T^{n+1}$ and $J^{n+1}$ from \eqref{semi-cor0};
    \\ \textrm{\bf Projection}: obtain $U^{n+1}$, $I^{n+1}$ from \eqref{proj00}.

\end{algorithm}

\begin{color}{black}
\subsection{Conservation of the time discretization}
Now we check the conservation of our three stage solver. Taking the integral with respect to $v$ in the first equation in \eqref{radiative_transport} and summing it up with the second equation, one gets the energy balance equation such that
\beq\label{energybalance}
\frac{1}{c}\partial_t \average{ I}+C_v\partial_t T+\frac{1}{\varepsilon } \average{v\cdot\nabla_x I}=0.
\eeq
A numerical discretization is called conservative if it can recover a discrete version of the energy balance equation. We show in the sequel that the semi-discretization proposed in this section is conservative.
From the correction step \eqref{semi-TI}, one gets
$$
 \f{1}{c}\left(ac\f{(T^{n+1})^4-(T^n)^4}{\Delta t}+\varepsilon\f{\f{C_v\eps}{\sigma\Delta t}(T^{n+1}-T^n)-\rho_J^n}{\Delta t}\right)+ C_v\f{T^{n+1}-T^n}{\Delta t}+ \f{1}{\varepsilon} \average{ v\cdot\nabla_x I^{n+1\ast} }=0.
$$
Subtracting the above equation by the integration of \eqref{semi-J} in $v$ and using \eqref{proj00} yields
$$
\f{\eps}{c}\f{\f{C_v\eps}{\sigma\Delta t}(T^{n+1}-T^n)-\rho_J^{n+1}}{\Delta t}+C_v\f{T^{n+1}-T^n}{\Delta t}+\f{1}{\eps}\average{v\cdot\nabla_x(I^{n+1\ast}-I^{n+1}) }=\f{\sigma}{\eps}\rho_J^{n+1}.
$$
 Therefore,
 $$
\big(1+\f{\eps^2}{c\sigma\Delta t}\big) \left(C_v\f{T^{n+1}-T^n}{\Delta t}-\f{\sigma}{\eps}\rho_J^{n+1}\right)+\f{1}{\eps}\average{v\cdot\nabla_x(I^{n+1\ast}-I^{n+1}) }=0.
 $$ 
 Combining the above equation with \eqref{semi-J} and the projection step \eqref{proj00}, we can obtain the following semi-discretization for \eqref{energybalance} such that
 \beq\label{semiconservative}
 \begin{aligned}
 &\f{1}{c \Delta t} \Big(\average{I^{n+1}}- \average{ I^n}\Big)  +C_v\f{T^{n+1}-T^n}{\Delta t}
 \\&+\frac{1}{\varepsilon }\average{v\cdot\nabla_x I^{n+1\ast} } -\frac{\eps}{c\sigma\Delta t+\eps^2} \average{v\cdot\nabla_x (I^{n+1\ast}-I^{n+1}) }=0.
 \end{aligned}
 \eeq
 When the solution is smooth enough, we have $|I^{n+1\ast}-I^{n+1}|\sim O(\Delta t^2)$.
 \eqref{semiconservative} is a consistent and conservative discretization for \eqref{energybalance}.  
 	\end{color}
\begin{remark}
There is another way of discretizing \eqref{radiative_transport3} by keeping the first term in \eqref{radiative_transport} unexpanded and discretizing it using the update from the prediction stage: 
\begin{subequations} \label{semi-cor}
\begin{numcases}{}
 \f{1}{c}\f{I^{n+1\ast}-I^n}{\Delta t}-v\cdot \nabla_x \left(\f{ac}{\sigma}v\cdot\nabla_x U^{n+1\ast} \right)+\varepsilon v\cdot\nabla_x J^{n+1}=-\sigma J^{n+1} \,; \\
C_v\frac{T^{n+1}-T^n}{\Delta t}= \sigma\rho_J^{n+1}. 
\end{numcases}
\end{subequations}
which avoids any nonlinear solver but yields a diffusion limit of the form 
\[
a\f{U^{n+1\ast}-U^n}{\Delta t}+C_v\f{T^{n+1}-T^n}{\Delta t}=\nabla_x \cdot \left(\f{ac}{3\sigma}\nabla_x U^{n+1\ast}\right)\,.
\]
We have tested the performance of this discretization, and it behaves well when $\sigma$ does not depend on T. However, it failed when $\sigma$ has a $T$ dependence. 
Besides,  \eqref{semi-cor} is {\it not} energy conservative.
Therefore, to facilitate the extension to the strongly nonlinear case such as $\sigma_T=\frac{\sigma}{T^3}$, we always use \eqref{radiative_transport3} from here on.

\end{remark}

\section{Spatial Discretization}
For the ease of exposition, we will explain our spatial discretion in 1D. That is, $x \in [0,L]$, $v \in [-1,1]$, and $\average{f(v)} = \frac{1}{2} \int_{-1}^1 f(v) \rd v$. The boundary condition \eqref{BC00} simplifies to 
\begin{equation} \label{BC1d}
I(t,0,v) = b_\text{L}(t,v), ~ \text{ for } ~v>0; \qquad I(t,L,v) = b_\text{R}(t,v), ~ \text{ for } ~ v<0\,.
\end{equation}
\black{The spatial domain is discretized as $-1= x_1 < x_2 < \cdots  < x_{N_x+1} = 1$. }
The higher dimensions can be treated in the dimension by dimension manner. To get a more compact stencil in spatial discretization, we use the even-odd parity method in the sequel.
\subsection{The prediction step}
Let
\begin{equation} \label{EO}
E(v)=\f{1}{2}\big(I(v)+I(-v)\big),  \qquad O(v)=\frac{1}{2\eps} (I(v)-I(-v)),\qquad  v>0
\end{equation}
be the even and odd part of $I$, then the time discretization of the prediction step in \eqref{radiative_IU} can be rewritten as
\begin{equation*}
\left\{
\begin{aligned}
& \frac{1}{c} \frac{O^{n+1\ast}-O^n}{\Delta t}+\frac{1}{\varepsilon^2} v\cdot \nabla_x E^{n+1\ast}=-\f{\sigma}{\varepsilon^2}O^{n+1\ast} \,;\\
&\frac{1}{c} \frac{E^{n+1\ast}-E^n}{\Delta t}+ v\cdot \nabla_x O^{n+1\ast}=\f{\sigma}{\varepsilon^2}\big(acU^{n+1\ast}-E^{n+1\ast}\big) \,;\\
&  C_v\frac{ U^{n+1\ast}-U^n}{\Delta t}=4\frac{\sigma}{\varepsilon^2}(T^n)^3\Big(\int_0^1E^{n+1\ast}\dv -acU^{n+1\ast}\Big).
\end{aligned}
\right.
\label{radiative_EO}
\end{equation*}
Consider the even part on half spatial grid and odd part on regular grid, i.e., 
\black{
\begin{align}
& E_{i+1/2} (v)\approx E(x_{i+1/2}, v), \quad i = 1, \cdots, N_x\,; \label{E000}
\\ & O_i (v)\approx O(x_i, v), \quad i = 1, \cdots, N_x + 1\,, \label{O000}
\end{align}
}
and we propose the following spatial discretization
\begin{subequations} \label{dis-pre}
\begin{numcases}{}
 \frac{1}{c} \frac{O_i^{n+1\ast}-O_i^n}{\Delta t}+\frac{1}{\varepsilon^2} v\frac{E^{n+1\ast}_{i+1/2}-E^{n+1\ast}_{i-1/2}}{\Delta x}=-\f{\sigma_i}{\varepsilon^2}O_i^{n+1\ast}, \quad  2 \leq i \leq N_x \,; \label{dis-pre1}\\
\frac{1}{c} \frac{E^{n+1\ast}_{i+1/2}-E^n_{i+1/2}}{\Delta t}+ v\frac{O^{n+1\ast}_{i+1}-O^{n+1\ast}_i}{\Delta x}=\f{\sigma_{i+1/2}}{\varepsilon^2}\big(acU_{i+1/2}^{n+1\ast}-E_{i+1/2}^{n+1\ast}\big), \quad 1 \leq i \leq N_x\,;\label{dis-pre2}\\
  C_v\frac{ U_{i+1/2}^{n+1\ast}-U_{i+1/2}^n}{\Delta t}=4\frac{\sigma_{i+1/2}}{\varepsilon^2}(T_{i+1/2}^n)^3\Big(\int_0^1E_{i+1/2}^{n+1\ast}\dv -acU_{i+1/2}^{n+1\ast}\Big)\,,  \quad  1 \leq i \leq N_x\,. \label{dis-pre3}
\end{numcases}
\label{radiative_EO_dis}
\end{subequations}
\black{To implement, we first express $U_{i+1/2}^{n+1\ast}$ in terms of $E_{i+1/2}^{n+1\ast}$ from \eqref{dis-pre3} and plug it into \eqref{dis-pre2}. Then \eqref{dis-pre2}, along with \eqref{dis-pre1}, forms a linear system for $E_{i+1/2}^{n+1\ast}$ and $O_{i+1/2}^{n+1\ast}$. Solving this system first and use the previous expression to get $U_{i+1/2}^{n+1\ast}$.}
\begin{remark}
Throughout the paper, we do not delineate on velocity discretization. In the numerical examples, we always use uniform grid in velocity with small enough grid size. 
\end{remark}

\subsection*{Boundary conditions} Note from \eqref{dis-pre1}--\eqref{dis-pre2} that one just need boundary conditions for $O_1^{n+1*}$ and $O_{N_x +1}^{n+1*}$. To cope with \eqref{BC1d}, using the relation \eqref{EO}, we have 
\begin{equation} 
E_{3/2}^{n+1*} + \frac{\eps}{2} \left( O_2^{n+1*}  + O_1^{n+1*} \right) = b_\text{L}(v), ~ v>0; \quad
E_{N_x+1/2}^{n+1*} - \frac{\eps}{2} \left( O_{N_x}^{n+1*}  + O_{N_x+1}^{n+1*} \right) = b_\text{R}(-v), ~ v>0\,.
\end{equation}
Hence,
\begin{equation} \label{bc-pre}
O_1^{n+1*}(v) = \frac{2}{\eps} \left( b_\text{L} (v) - E_{3/2} ^{n+1*}\right) - O_2^{n+1*}, \quad 
O_{N_x+1}^{n+1*}(v) = \frac{2}{\eps} \left( E_{N_x + \half} ^{n+1*}-b_\text{R} (-v) \right) - O_{N_x}^{n+1*}
\end{equation}

\subsection*{Diffusion limit of \eqref{radiative_EO_dis}:} First from \eqref{dis-pre1} and \eqref{dis-pre2}, one sees that
\begin{equation} \label{9162}
O_i^{n+1*} =  - \frac{v}{\sigma_i \Dx } \left( E_{i + \half}^{n+1*} - E_{i-\half}^{n+1*}\right) + \order (\eps^2), \qquad 
E_{i+\half}^{n+1*} = ac U_{i+\half}^{n+1*} + \order(\eps^2)\,.
\end{equation}
Taking average in $v$ of \eqref{dis-pre2}, then adding it to \eqref{dis-pre3} divided by $4(T_{i+\half}^{n})^3$, we have
\[
\frac{1}{c \Dt} \left(\average{E_{i+\half}^{n+1*} }- \average{E_{i+\half}^n} \right) + \frac{C_v}{4\Dt (T_{i+\half}^n)^3} \left( U_\ihalf^{n+1*} - U_\ihalf^n \right)  + \average{ \frac{v}{\Dx} \left( O_{i+1}^{n+1*} - O_i^{n+1*}\right)} = 0 \,.
\]
Plugging the relation in \eqref{9162}, it becomes
\[
\left( \frac{a}{\Dt} + \frac{C_v}{4\Dt (T_\ihalf^n)^3}\right) \left(  \Uns_\ihalf - U_{i+\half}^n \right) 
+ \average{\frac{ac v}{\Dx} \left[  -\frac{v}{\sigma_{i+1} \Dx} (U_{i+\frac{3}{2}}^{n+1*} - U_\ihalf^{n+1*}) 
                                               + \frac{v}{\sigma_i \Dx} (U_\ihalf^{n+1*} - U_{i-\half}^{n+1*})\right]} = 0\,
\]
to the leading order, which is the same as the first equation in \eqref{dis-diff}.

{\color{black}\subsection{The correction step} In the correction step, we always have $U=T^4$.
\black{Plugging the expansion \eqref{J00} into \eqref{EO}, we have 
\[
E = acT^4 + \f{1}{2}\big(J(v)+J(-v)\big) := acT^4 + \eps E_J, \quad O = \f{1}{2\eps}(\eps J(v) - \eps J(-v)) :=O_J\,.
\]
Consequently, \eqref{radiative_transport3} can be written into}
\begin{equation} \label{9163}
\left\{
\begin{aligned}
& \frac{1}{c} \partial_t O  + \f{v}{\eps^2} \cdot \nabla_x \left( acU \right)+ \f{v}{\eps}\cdot  \nabla_x E_J = -\frac{\sigma}{\eps^2} O_J\,;
\\ & \frac{1}{c} \partial_t E+ v \cdot \nabla_x O_J = - \f{\sigma}{\eps} E_J\,;
\\ &  C_v \partial_t T = \f{\sigma}{\eps} \rho_J\,.
\end{aligned}
\right.
\end{equation}
and the semi-discretization \eqref{semi-cor10} becomes
\begin{subequations} \label{dis-semi-cor}
\begin{numcases}{}
 \frac{1}{c} \frac{\tilde O^{n+1} - O^n}{\Delta t}+ \f{v}{\eps^2} \cdot \nabla_x \left( acU^{n+1\ast} \right)+\f{v}{\eps}\cdot\nabla_x E_J^{n+1\ast}=-\frac{\sigma}{\varepsilon^2} \tilde O^{n+1}_{J} \,;\nonumber\\
\frac{1}{c} \frac{\tilde E^{n+1} - E^n}{\Delta t}+ v\cdot\nabla_x O_J^{n+1\ast} =-\f{\sigma}{\eps}  \tilde E^{n+1}_{J}\,; \nonumber\\
 C_v\frac{ T^{n+1}-T^n}{\Delta t}=\f{\sigma}{\eps} \int_0^1E^{n+1}_{J}\dv\,. \nonumber
\end{numcases}
\end{subequations}
The spatial discretization then takes the following form: 
\begin{subequations} \label{dis-cor}
\begin{numcases}{}
 \frac{1}{c} \frac{\tilde O_i^{n+1}-O_i^n}{\Delta t}+\frac{acv}{\eps^2\Delta x}(U^{n+1\ast}_{i+1/2}-U^{n+1\ast}_{i-1/2})+\frac{v}{\eps}\frac { E^{n+1\ast}_{Ji+1/2}- E^{n+1\ast}_{Ji-1/2}}{\Delta x}\nonumber \\
\hspace{7cm} =-\frac{\sigma_i}{\varepsilon^2} \tilde O^{n+1}_{Ji}, \quad 2 \leq i \leq N_x \,; \label{dis-cor1}\\
\frac{1}{c} \frac{\tilde E^{n+1}_{i+1/2}-E^n_{i+1/2}}{\Delta t} + v\frac{O^{n+1*}_{Ji+1}-O^{n+1*}_{Ji}}{\Delta x}  =-\f{\sigma_{i+1/2}}{\eps} \tilde E^{n+1}_{Ji+1/2}, \quad 2 \leq i \leq N_x-1\,; \label{dis-cor2}\\
 C_v\frac{ T_{i+1/2}^{n+1}-T_{i+1/2}^n}{\Delta t}=\f{\sigma_{i+1/2}}{\eps}\int_0^1\tilde E^{n+1}_{Ji+1/2}\dv\,, \quad 1\leq i \leq N_x \,. \label{dis-cor3}
\end{numcases}
\end{subequations}
where the $\tilde E_\ihalf^{n+1}$ and $\tilde O_i^{n+1}$ appeared above are replaced by
\begin{equation} \label{dis-cor4}
\tilde E_{i+\half}^{n+1} = ac (T_{i+\half}^{n+1})^4 + \eps \tilde E_{J\ihalf}^{n+1}\,, \quad 
\tilde O_i^{n+1} =  \tilde{O}_{Ji}^{n+1}\,.
\end{equation}

\black{To implement \eqref{dis-cor}, we first integrate \eqref{dis-cor2} with respect to $v$ and add it to \eqref{dis-cor3}, so that the right hand side vanishes. From \eqref{dis-cor4} and since $\int \tilde E_{J i +\half}^{n+1} \rd v$ can be written in terms of $T_{i+\half}^{n+1}$ from \eqref{dis-cor3}, it leads to an equation only for $T_{i+\half}^{n+1}$. That is,\begin{align}
& \frac{1}{c\Delta t} \left[ac(T_{i+\half}^{n+1})^4 + \eps^2 \frac{C_v}{\sigma_{i+\half}\Delta t} (T_{i+\half}^{n+1} - T_{i+\half}^{n}) \right] 
+ \frac{C_v}{\Delta t} \left(T_{i+\half}^{n+1} - T_{i+\half}^{n} \right) \nonumber
\\ & \hspace{7cm} =  \int_0^1 \left[  \frac{1}{c\Delta t}E^n_{i+1/2} -  v\frac{O^{n+1 *}_{Ji+1}-O^{n+1*}_{Ji}}{\Delta x}  \right] \rd v\,.  \label{Tnp1}
\end{align}
 Here $O_{Ji}^{n+1*} = O_i^{n+1*}$ and is known from the prediction step.
 }

\black{
Once $T_{i+1/2}^{n+1}$ is obtained, we replace $U^{n+1\ast}_{i\pm 1/2}$ by $(T^{n+1}_{i\pm 1/2})^4$; $O_J^{n+1\ast}$, $\tilde O_J^{n+1}$ by $O_J^{n+1}$ and $E_J^{n+1\ast}$, $\tilde E_J^{n+1}$ by $E_J^{n+1}$ in \eqref{dis-cor1} and \eqref{dis-cor2} and consider it as a linear system for $O_{Ji}^{n+1}$ and $E_{Ji+\half}^{n+1}$.  
}

\subsection*{Boundary conditions}\black{ Since \eqref{Tnp1} does not need a boundary condition, we only need the boundary conditions in \eqref{dis-cor1}--\eqref{dis-cor2} for $O_1^{n+1}$ and $O_{N_x +1}^{n+1}$. As with \eqref{bc-pre}, we have 
\[
O_{J1}^{n+1}=O_1^{n+1}(v) = \frac{2}{\eps} \left( b_\text{L} (v) - E_{3/2} ^{n+1}\right) - O_2^{n+1}, \quad 
O_{JN_x+1}^{n+1} =O_{N_x+1}^{n+1}(v) = \frac{2}{\eps} \left( E_{N_x + \half} ^{n+1}-b_\text{R} (-v) \right) - O_{N_x}^{n+1}\,.
\]
Using these boundary conditions, we have additional equations to \eqref{dis-cor2} for $i=1$ and $N_x$ to update $E_{3/2}^{n+1}$ and $E_{N_x+1/2}^{n+1}$.
}

\subsection*{Diffusion limit of \eqref{dis-cor}} \black{Note from \eqref{dis-pre1} and \eqref{dis-pre2} , one can see that
\[
O_{J i} ^{n+1*}=O_{ i} ^{n+1*}=-\f{acv}{\sigma_i}\f{\Uns_\ihalf - \Uns_{i-\half}}{\Delta x}+\order(\eps^2).
\]}
 Taking the average in $v$ of \eqref{dis-cor2} and adding it to \eqref{dis-cor3}, keeping the leading order terms, we have
\[
\frac{1}{c}\frac{\average{E_\ihalf^{n+1} - E_\ihalf^n}}{\Dt} + C_v \frac{T_\ihalf^{n+1} - T_\ihalf^n}{\Dt} - \frac{ac}{3 \Dx^2} \left( \frac{\Uns_{i+\frac{3}{2}} - \Uns_\ihalf}{\sigma_{i+1}}  - \frac{\Uns_\ihalf - \Uns_{i-\half}}{\sigma_i}  \right) = 0\,.
\]
Then using the relation between $E^{n+1}_{\ihalf}$ and $U^{n+1}_\ihalf$ from \eqref{dis-cor4}, the above formula is the same as the second step in \eqref{dis-diff}. 
}

\subsection{The projection step}
$U$, $O$ and $E$ for the next time step are determined by
\begin{subequations} \label{dis-proj}
\begin{numcases}{}
U_{i+1/2}^{n+1}=(T_{i+1/2}^{n+1})^4, \qquad 1 \leq i \leq N_x\,; \label{dis-proj1}
\\ \black{O^{n+1}_i= O^{n+1}_{Ji},} \qquad 2 \leq i \leq N_x \,;  \label{dis-proj2}
\\ \black{E^{n+1}_{i+1/2}=acU^{n+1}_{i+1/2}+\varepsilon E^{n+1}_{Ji+1/2},} \qquad 1 \leq i \leq N_x\,.  \label{dis-proj3}
\end{numcases}
\end{subequations}
Here the boundary condition for $O^{n+1}$ in \eqref{dis-proj2} is similar to that in \eqref{bc-pre}, i.e., 
\begin{equation*} 
O_1^{n+1}(v) = \frac{2}{\eps} \left( b_\text{L} (v) - E_{3/2} ^{n+1}\right) - O_2^{n+1}, \quad 
O_{N_x+1}^{n+1}(v) = \frac{2}{\eps} \left( b_\text{R} (-v) - E_{N_x + \half} ^{n+1}\right) - O_{N_x}^{n+1}\,.
\end{equation*}
To summarize, we have the following one time step update of the fully discrete version of GRTE.
 \begin{algorithm}[!h]
\caption{one step of fully discrete update for GRTE}
\SetAlgoLined
\KwIn{$E_{i+\half}^n$, $O^n_i$, $U_{i+\half}^n$, $T_{i+\half}^n$, ($U_{i+\half}^n = (T_{i+\half}^{n})^4$) }
\KwOut{$E_{i+\half}^{n+1}$, $O^{n+1}_i$, $U_{i+\half}^{n+1}$, $T_{i+\half}^{n+1}$, ($U_{i+\half}^{n+1} = (T_{i+\half}^{n+1})^4$) }

\BlankLine
    \textrm{\bf Prediction}: obtain $O_i^{n+1*}$, $E_{i+\half}^{n+1*}$, $U_{i+\half}^{n+1*}$ from \eqref{dis-pre}; 
    \\ \textrm{\bf Correction}: obtain $O_{Ji}^{n+1}$, $E_{J i+\half}^{n+1}$, $T_{i+\half}^{n+1}$ from \eqref{dis-cor}:
    \\ \hspace{1.3cm} \black{\underline{substep1}:  get $T_{i+\half}^{n+1}$ from \eqref{Tnp1}};
    \\  \hspace{1.3cm} \black{\underline{substep2}:  get $E_{J i+\half}^{n+1}$ and $O_{J i+\half}^{n+1}$ from \eqref{dis-cor1} and \eqref{dis-cor2} with \eqref{dis-cor4}, and $U^{n+1\ast}_i\pm 1/2$ replaced by $(T^{n+1}_{i\pm 1/2})^4$; $O_J^{n+1\ast}$, $\tilde O_J^{n+1}$ replace by $O_J^{n+1}$ and $E_J^{n+1\ast}$, $\tilde E_J^{n+1}$ replace by $E_J^{n+1}$ };
    \\ \textrm{\bf Projection}: obtain $U_{i+\half}^{n+1}$, $O_i^{n+1}$ and $E_{i+\half}^{n+1}$ from \eqref{dis-proj}.
\end{algorithm}

\black{We now calculate the computational complexity of our method. Assume that in velocity direction we have $N_v$ cells, and we compute the solution up to time $N_t \Delta t$. Then from \eqref{E000} and \eqref{O000} we see that $E$ and $O$ are vectors of size $N_x N_v$ and $(N_x+1)N_v$ respectively. The following computations are needed for each time step:
\begin{itemize}
\item Prediction: solving $(2N_x+1)N_v$ linear system to get $E^{n+1*}$ and $O^{n+1*}$, and then $N_x$ explicit scalar updates to get $U^{n+1*}$;
\item Correction: $N_x$ scalar Newton iterations to get $T^{n+1}$, and solving a $(2N_x+1)N_v$ linear system to get $E_J^{n+1}$ and $O_J^{n+1}$;
\item Projection: $(2N_x +1) N_v$ explicit scalar updates to get $O^{n+1}$ and $E^{n+1}$.
\end{itemize}}

\section{\black{The nonlinearity in the opacity}}
As with the previous three-stage scheme, we describe below the corresponding stages for the \black{nonlinear opacity} \eqref{marshak00} with emphasis on the variations compared to the previous case. The method we developed here can be used for other kinds of nonlinearity, such as non-grey radiative transfer equation.

\subsection{The prediction step}
Noting that $\sigma_T=\sigma/T^3$, similar to the prediction stage \eqref{pre00}, we solve the following system 
\begin{equation} \label{pre-marshak}
\left\{\begin{aligned}
&\frac{1}{c} \partial_t I+\frac{1}{\varepsilon} v\cdot \nabla_x I=
\frac{1}{\varepsilon^2} \frac{\sigma}{T^3}\left(acU-I\right)\,; \\
&C_v\partial_t U=4\frac{\sigma}{\varepsilon^2}(\rho-acU).
\end{aligned} \right.
\end{equation}
Let $E$ and $O$ be defined the same as in \eqref{EO}, and denote 
$$K=T^3\,,$$
we discretize \eqref{pre-marshak} as 
\begin{subequations} \label{dis-pre-marshak}
\begin{numcases}{}
K_i^n\frac{1}{c} \frac{O_i^{n+1\ast}-O_i^n}{\Delta t}+K_i^n\frac{1}{\varepsilon^2} v\frac{E^{n+1\ast}_{i+1/2}-E^{n+1\ast}_{i-1/2}}{\Delta x}=-\f{\sigma_{i}}{\varepsilon^2}O_i^{n+1\ast}, \quad 2 \leq i \leq N_x\,;  \label{dis-pre-marshak1}\\
K_{i+1/2}^n\frac{1}{c} \frac{E^{n+1\ast}_{i+1/2}-E^n_{i+1/2}}{\Delta t}+ K_{i+1/2}^n v\frac{O^{n+1\ast}_{i+1}-O^{n+1\ast}_i}{\Delta x}
=\f{\sigma_{i+1/2}}{\varepsilon^2}\big(acU_{i+1/2}^{n+1\ast}-E_{i+1/2}^{n+1\ast}\big), ~  1\leq i \leq N_x\,; \qquad \\
  C_v\frac{ U_{i+1/2}^{n+1\ast}-U_{i+1/2}^n}{\Delta t}=4\frac{\sigma_{i+1/2}}{\varepsilon^2}\Big(\int_0^1E_{i+1/2}^{n+1\ast}\dv -acU_{i+1/2}^{n+1\ast}\Big), \quad 1 \leq i \leq N_x\,,
\end{numcases}
\end{subequations}
where $V_{i+1/2}^n$, $V_i^n$ are respectively the approximations to $(T^n(x_{i+1/2}))^{3}$ and $(T^n(x_i))^{3}$, which are determined by 
\[
K_{i+1/2}^n=(T_{i+1/2}^n)^3,\qquad K_i=\f{1}{2}(K_{i+1/2}+K_{i-1/2}).
\]
Compared to \eqref{dis-pre}, the major difference here is the presence of $K$, and we prefer to multiply it to the left hand side just to avoid dividing by zero when $K=0$. The boundary conditions $O_1^{\nps}$ and $O_{N_x+1}^\nps$ needed in \eqref{dis-pre-marshak1} is chosen the same as in \eqref{bc-pre}.

\subsection{The correction step}
First we derive the correction system. As with \eqref{J00}, we expand
\begin{color}{black}{
\[
I = acU + \eps J\,.
\]
 We define 
\begin{equation*} 
E(v)=acU+\varepsilon E_J(v), \qquad O(v)= O_J(v)\,.
\end{equation*}
Then we get the following correction system similar to \eqref{9163} as

\begin{equation*} 
\left\{
\begin{aligned}
& \frac{K}{c} \partial_t O + \f{K}{\eps^2} v \cdot  \nabla_x(acU)+ \f{K}{\eps}v \cdot  \nabla_x E_J = -\frac{\sigma}{\eps^2} O_J\,,
\\ & \frac{K}{c} \partial_t E  +  K v \cdot \nabla_x O_J = - \f{\sigma}{\eps} E_J\,,
\\ & K C_v \partial_t T = \f{\sigma}{\eps}\rho_J\,,
\end{aligned}
\right.
\end{equation*}
which can be discretized as follows 
\begin{subequations} \label{dis-cor-marshak}
\begin{numcases}{}
 \frac{\Vis}{c} \frac{\tilde O_i^{n+1}-O_i^n}{\Delta t}+\f{acv\Vis}{\eps^2} \frac {U^{n+1*}_{i+1/2}-U^{n+1*}_{i-1/2}}{\Delta x}+\f{\Vis}{\eps}v\frac {E^{n+1\ast}_{Ji+1/2}-E^{n+1\ast}_{Ji-1/2}}{\Delta x}\nonumber \\
\hspace{10cm} =-\frac{\sigma_{i}}{\varepsilon^2} \tilde O^{n+1}_{Ji}, \quad 2 \leq i \leq N_x \,; \label{dis-cor-marshak1}\\
\frac{\Vips}{c} \frac{\tilde E^{n+1}_{i+1/2}-E^n_{i+1/2}}{\Delta t}+ \Vips  v\frac{O_{Ji+1}^{\color{black}n+1*}-O_{Ji}^{\color{black}n+1*}}{\Delta x} =-\f{\sigma_{i+1/2}}{\eps} \tilde E^{n+1}_{Ji+1/2}, \quad 2 \leq i \leq N_x-1\,; \label{dis-cor-marshak2}\\
\Vips C_v\frac{ T_{i+1/2}^{n+1}-T_{i+1/2}^n}{\Delta t}=\f{\sigma_{i+1/2}}{\eps}\int_0^1\tilde E^{n+1}_{Ji+1/2}\text{d}v\,, \quad 1\leq i \leq N_x \,. \label{dis-cor-marshak3}
\end{numcases}
\end{subequations}
Like before, $\tilde E_\ihalf^{n+1}$ and $\tilde O_i^{n+1}$ are replaced by
\begin{equation}  \label{9171}
\tilde E_{i+\half}^{n+1} = ac (T_{i+\half}^{n+1})^4 + \eps \tilde {E}_{J\ihalf}^{n+1}\,, \quad 
\tilde O_i^{n+1} =  \tilde{O}_{Ji}^{n+1}\,, 
\end{equation}
Here $\Vips = (U_{i+\half}^{n+1*} )^{3/4}$, and $\Vis = \frac{1}{2} (\Vips + V_{i-\half}^{n+1*})$. \black{It is important to notice that as discussed in section 3.2, $T$ is updated first by scalar Newton iterations, while in (\ref{dis-cor-marshak}), at the grid points such that $\Vips=0$, $T_{i+\frac{1}{2}}^{n+1}$ can not be determined. Therefore, to update $T$, we have to consider}
\begin{subequations}
\begin{align}
\frac{1}{c} \frac{E^{n+1}_{i+1/2}-E^n_{i+1/2}}{\Delta t}&  +  v\frac{O^{n+1*}_{Ji+1}-O^{n+1*}_{Ji}}{\Delta x} =-\frac{\sigma_{i+1/2}}{\eps\Vips } \tilde E^{n+1}_{Ji+1/2}, \quad 2 \leq i \leq N_x-1\,;\label{dis-cor-marshak4}\\
&C_v\frac{ T_{i+1/2}^{n+1}-T_{i+1/2}^n}{\Delta t}=\frac{\sigma_{i+1/2}}{\eps\Vips }\int_0^1\tilde E^{n+1}_{Ji+1/2}\text{d}v\,, \quad 1\leq i \leq N_x \,. \label{dis-cor-marshak5}
\end{align}
\end{subequations}
Though the right hand sides of the above equations are meaningless when $\Vips=0$, integrate (\ref{dis-cor-marshak4}) with respect to $v$, summing it up with (\ref{dis-cor-marshak5}), the right hand sides canceled with each other and we obtain the equation to get $T^{n+1}$. After getting $T^{n+1}$, one can replace 
$K^{n+1\ast}$ by $(T^{n+1})^3$, $U^{n+1\ast}$ by $(T^{n+1})^4$, both $E_J^{n+1\ast}$ and $\tilde E_J^{n+1}$ by $E_J^{n+1}$; both $O_J^{n+1\ast}$ and $\tilde O_J^{n+1}$ by $O_J^{n+1}$ in \eqref{dis-cor-marshak1}-\eqref{dis-cor-marshak2}, then solve a linear system for $E_J^{n+1}$ and $O_J^{n+1}$.
}\end{color}

\subsection{The projection step}
Finally, $U$, $E$, $O$ for the next time step are determined by 
\begin{subequations}
\begin{numcases}{}
U_{i+1/2}^{n+1}=(T_{i+1/2}^{n+1})^4, \qquad 1 \leq i \leq N_x\,; \nonumber
\\ \black{O^{n+1}_i= O_{Ji}^{n+1}}, \qquad 2 \leq i \leq N_x \,;  \nonumber
\\ \black{E^{n+1}_{i+1/2}=acU^{n+1}_{i+1/2}+\varepsilon E^{n+1}_{Ji+1/2},} \qquad 1 \leq i \leq N_x\,.  \nonumber
\end{numcases}
\end{subequations}

\begin{remark}
The corresponding prediction-correction-projection for the diffusion limit takes the form
\begin{subequations}
\begin{numcases}{}
\frac{U^{n+1\ast}-U^n}{\Delta t}=\frac{4(T^n)^3}{4a(T^{n})^3+C_v} \nabla_x \cdot \vpran{\frac{ac (T^n)^3}{\sigma \dc}\nabla_x U^{n+1\ast}} \nonumber\,;
\\ a\frac{(T^{n+1})^4-U^n}{\Delta t}+C_v \frac{T^{n+1}-T^n}{\Delta t}=\nabla_x \cdot \vpran{\frac{ac (T^{n+1*})^3}{\sigma \dc}\nabla_x U^{n+1^\ast}},\quad U^{n+1} = (T^{n+1})^4\,,  \nonumber
\end{numcases}
\end{subequations}
which is consistent with \eqref{marshak-diff2}.
\end{remark}

\section{Numerical Examples}
In this section, we conduct a few numerical experiments to test the performance of our proposed method. Our examples will cover both optically thin $\eps \sim \mathcal{O}(1)$ and optically thick $\eps \ll 1$ cases, and with both strong $\sigma \sim \mathcal{O}(1)$ and weak $\sigma \ll 1$ scattering. Without loss of generality, we assume $a = c= C_v = 1$ throughout the examples except the last one. 
\subsection{Three-stage diffusion solver}
\begin{figure}[!h] 
\includegraphics[width=0.45\textwidth]{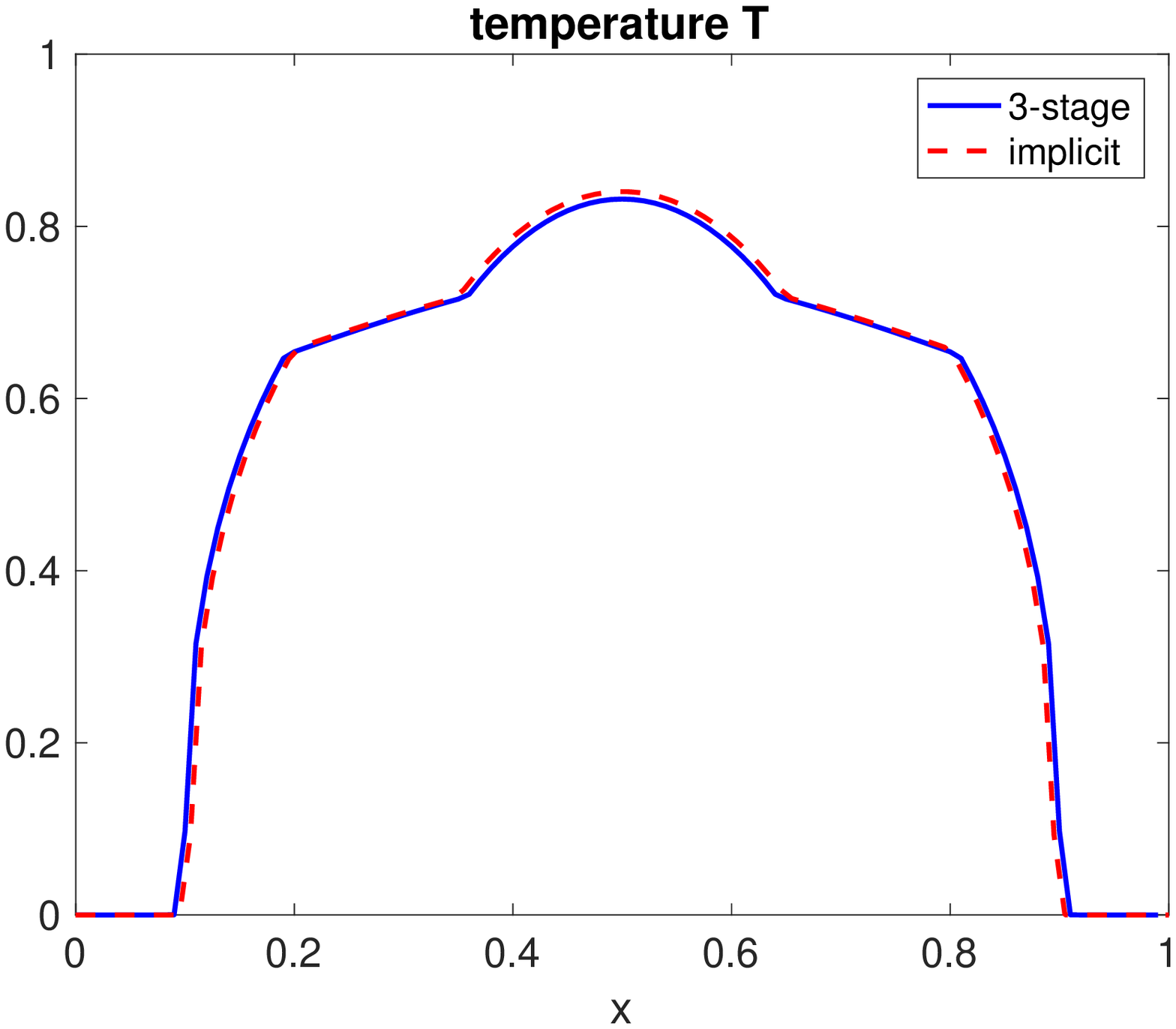}
\includegraphics[width=0.45\textwidth]{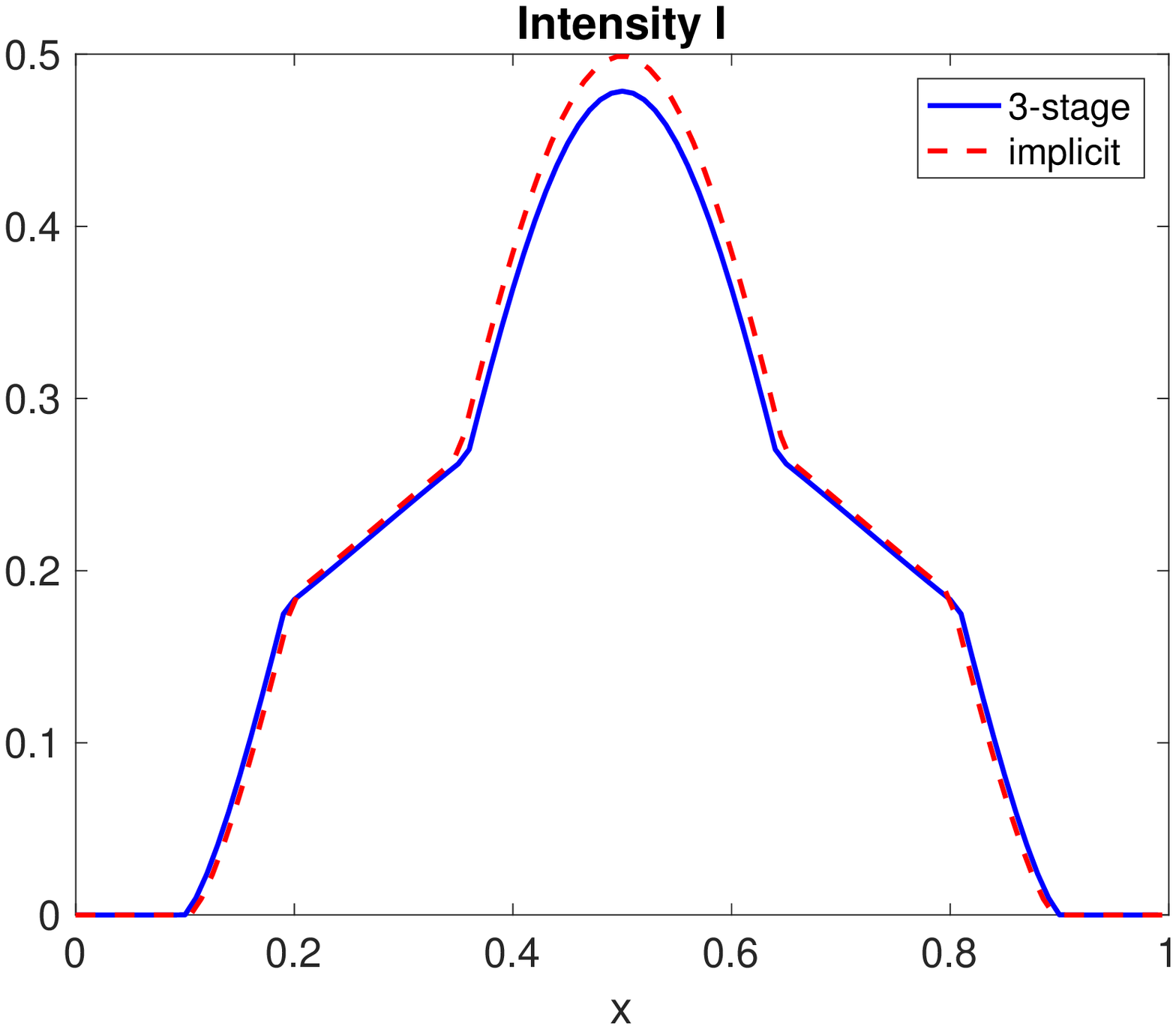}
\caption{For both methods---our new method \eqref{eq:limit_U_dis1}--\eqref{eq:limit_UT_dis3} and the fully implicit solver \eqref{imp_T}, $\Delta x = 0.01$, $\Delta t = 0.05\Delta x$. Left: temperature $T$. Right:  $ac T^4$. }
\label{fig-diff}
\end{figure}

The limit nonlinear diffusion equation is degenerate when $T=0$, thus there exist special non-negative compact supported solutions. In the literature, these kind of solutions serve as benchmark numerical tests. In this subsection, we first check the accuracy of our three-stage diffusion solver \eqref{eq:limit_U_dis1}--\eqref{eq:limit_UT_dis3}, which only needs scalar Newton's solver. As a reference, we use a fully implicit scheme for 
\eqref{eq:limit_T}:
\begin{align} \label{imp_T}
a \frac{(T_i^{n+1})^4 - (T_i^{n})^4}{\Delta t} + C_v \frac{T_i^{n+1}-T_i^n}{\Delta t} = \frac{ac}{3\Delta x^2} 
\left( \frac{(T_{i+1}^{n+1})^4 - (T_i^n)^4}{\sigma_{i+1/2}}  - \frac{(T_i^{n+1})^4 - (T_{i-1}^{n+1})^4}{\sigma_{i-1/2}} \right)\,.
\end{align}
The initial condition is chosen as
\[
T(x,0) = \max\{\sin(2\pi(x-1/4)),0 \}^{\frac{1}{4}}, \quad  x \in [0,1]\,,
\]
and the collision cross-section $\sigma$ takes the form
\begin{equation} \label{eqn:sigma}
\sigma(x) = \left\{ \begin{array}{cc} \sigma_0 & 0.2\leq x \leq 0.35, ~ 0.65\leq x \leq 0.8\\ 1 & \text{elsewhere} \end{array} \right. \,,
\end{equation}
where $\sigma_0 = 0.2$. The results are plotted in Fig.~\ref{fig-diff}, where good agreements between our 3-stage solver and fully implicit solver are observed. \black{We also note that when we further decrease $\Delta x$, we can keep the relation $\Delta t = 0.05 \Delta x$ unchanged, which confirms a hyperbolic CFL condition.}

\subsection{Accuracy}\label{section6.2}
\begin{figure}[!h]
\includegraphics[width=0.45\textwidth]{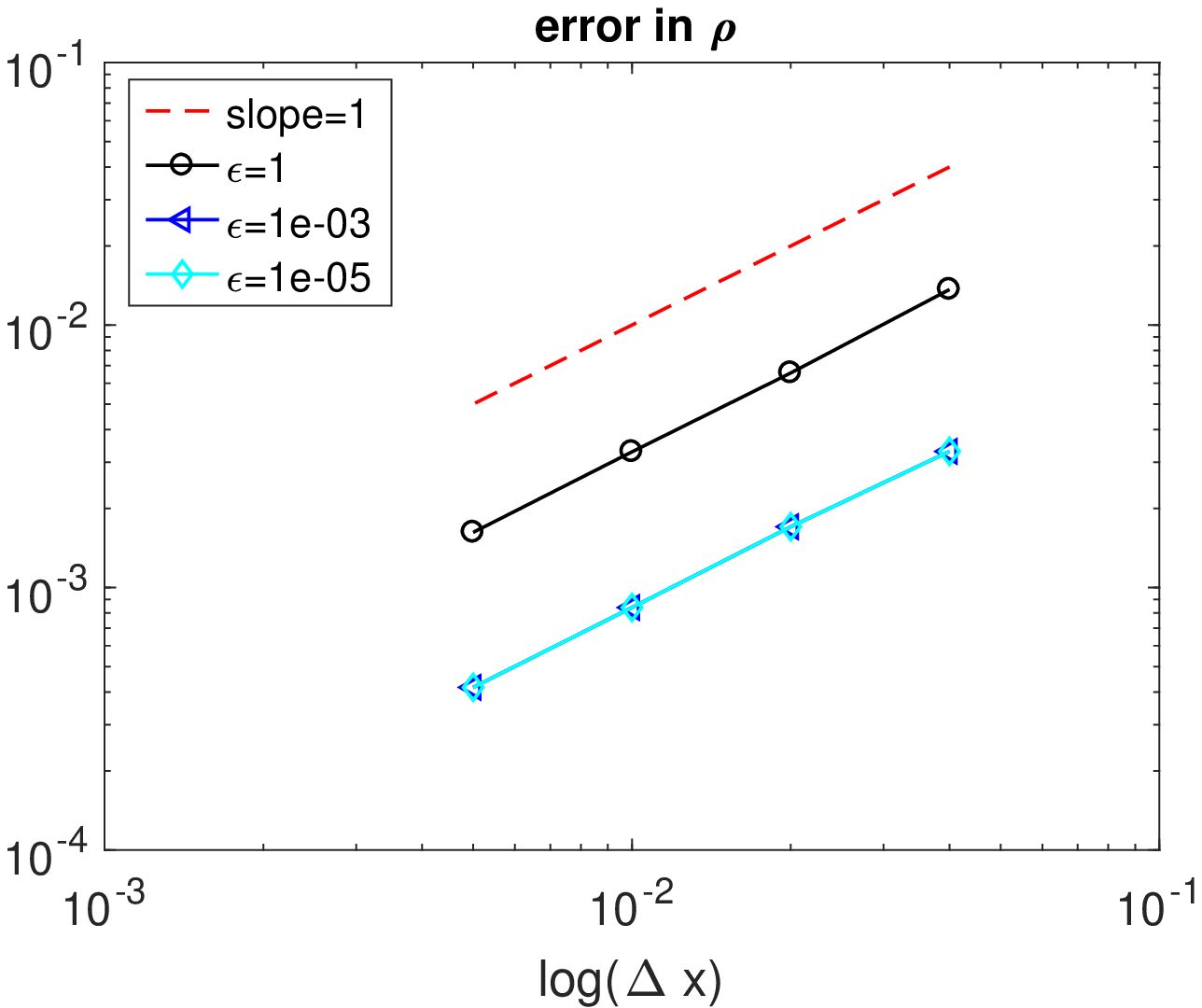}
\includegraphics[width=0.45\textwidth]{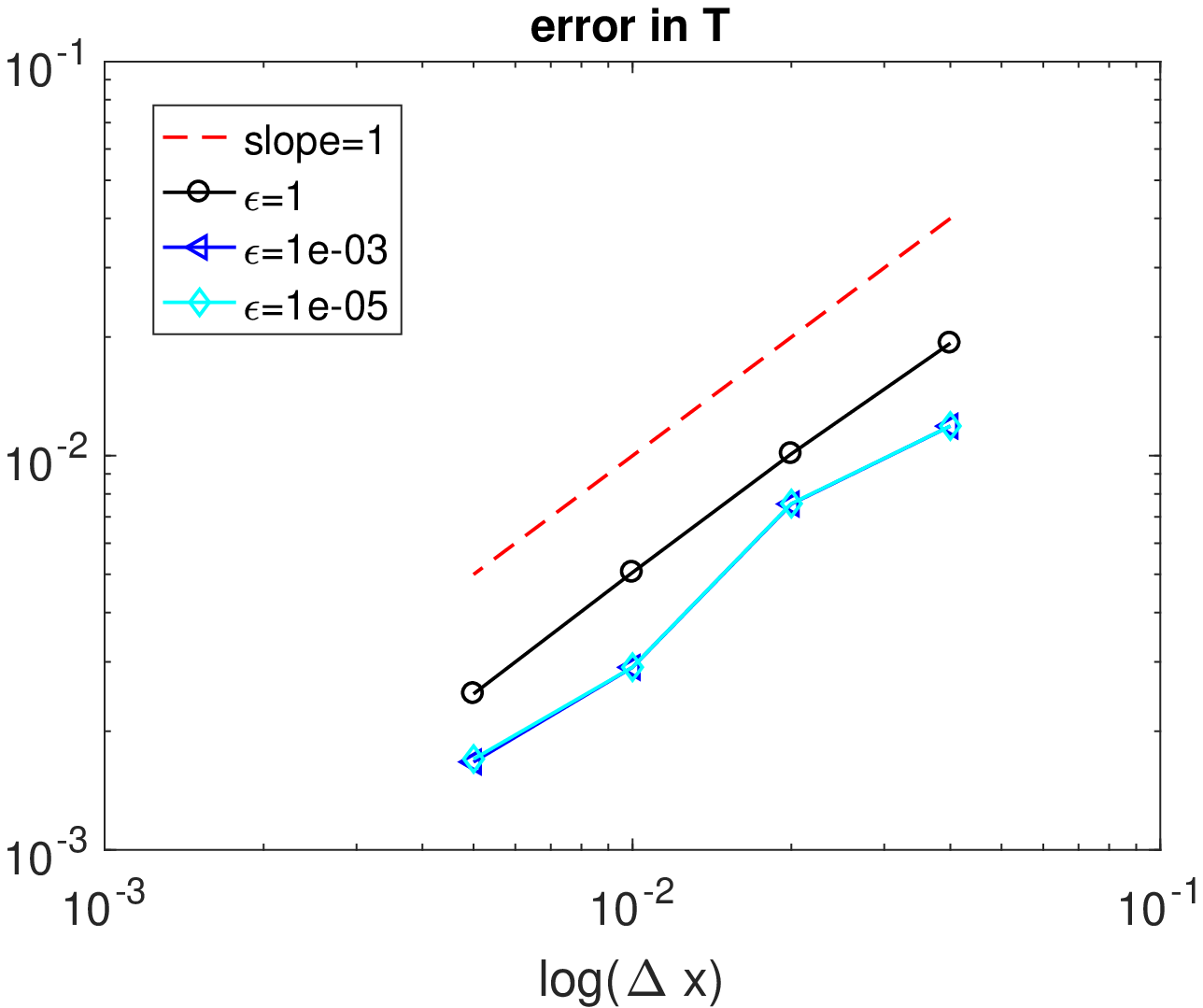}
\caption{Plot of error \eqref{error} with decreasing $\Delta x = \frac{1}{25}$, $\frac{1}{50}$, $\frac{1}{100}$, $\frac{1}{200}$,$\frac{1}{400}$. Here $\Delta t = 0.1\Delta x$.}
\label{fig:conv}
\end{figure}

This subsection aims at testing the accuracy of our AP scheme. The computational domain is $[0,1]$ and
the initial temperature and density fluxes are
\begin{equation} \label{IC000}
T(0,x)=\max\big\{1-40(x-1/2)^2,0\big\},\qquad
I(0,x, v)=acT(0,x)^4\,.
\end{equation}
The collision cross-section is chosen as $\sigma(x)=1$. In Figure \ref{fig:conv}, we plot error 
\begin{equation} \label{error}
\text{error}_\rho = \| \rho_{\Delta x}(\cdot, t_\text{max}) -\rho_{\Delta x/2}(\cdot, t_\text{max})\|_{l_1},\quad
\text{error}_T = \| T_{\Delta x}(\cdot, t_\text{max}) - T_{\Delta x/2}(\cdot, t_\text{max})\|_{l_1}
\end{equation}
with decreasing $\Delta x$ and $t_\text{max} = 0.1$, for different values of $\varepsilon = 1$, $10^{-3}$, $10^{-5}$. A uniform first order accuracy across different regimes is verified.  
\begin{color}{black}{
\subsection{Stability condition}
To check the the stability, we record the choice of $\Delta t$ for different opacity in Table \ref{tab:stability1} and Table \ref{tab:stability2}  with different $\eps$ and $\Delta x$ by using the same initial condition in Section \ref{section6.2}. Here we choose $\Delta t=C\Delta x$  with C being the largest constant such that the scheme is stable, and record $C$ in the table. 
\begin{table}    
 \centering
    \fontsize{10}{12}\selectfont    
    \caption{Stability test for $\sigma=1$.}
\begin{tabular}{l|ccccc}
    \hline
\diagbox[width=5em,trim=l] {$\eps$}{$\Delta x$} & $\f{1}{25}$ & $\f{1}{50}$ & $\f{1}{100}$  & $\f{1}{200}$ & $\f{1}{400}$ \\
\hline
1 & $>10$ & $>10$ & $>10$ & $>10$ & $>10$  \\
1e-03 & 0.8 & 0.5 & 0.5 & 0.5 & 0.5 \\
1e-05 & 0.8 & 0.5 & 0.5 & 0.5 & 0.5  \\
    \hline
\end{tabular}\vspace{0cm}
   \label{tab:stability1}
    \end{table}
 \begin{table}
    \centering
    \fontsize{10}{12}\selectfont    
    \caption{Stability test for $\sigma_T=\f{1}{T^3}$.}
\begin{tabular}{l|ccccc}
    \hline
\diagbox[width=5em,trim=l] {$\eps$}{$\Delta x$} & $\f{1}{25}$ & $\f{1}{50}$ & $\f{1}{100}$  & $\f{1}{200}$ & $\f{1}{400}$ \\
\hline
1 & 16 & 10 & 8 & 6 & 6  \\
1e-03 & 2 & 2 & 1 & 1 & 1 \\
1e-05 & 2 & 2 & 1 & 1 & 1 \\
    \hline
\end{tabular}\vspace{0cm}
    \label{tab:stability2}
\end{table}
}\end{color}
\subsection{Varying Sigma}
\begin{figure}[!h]
\includegraphics[width=0.45\textwidth]{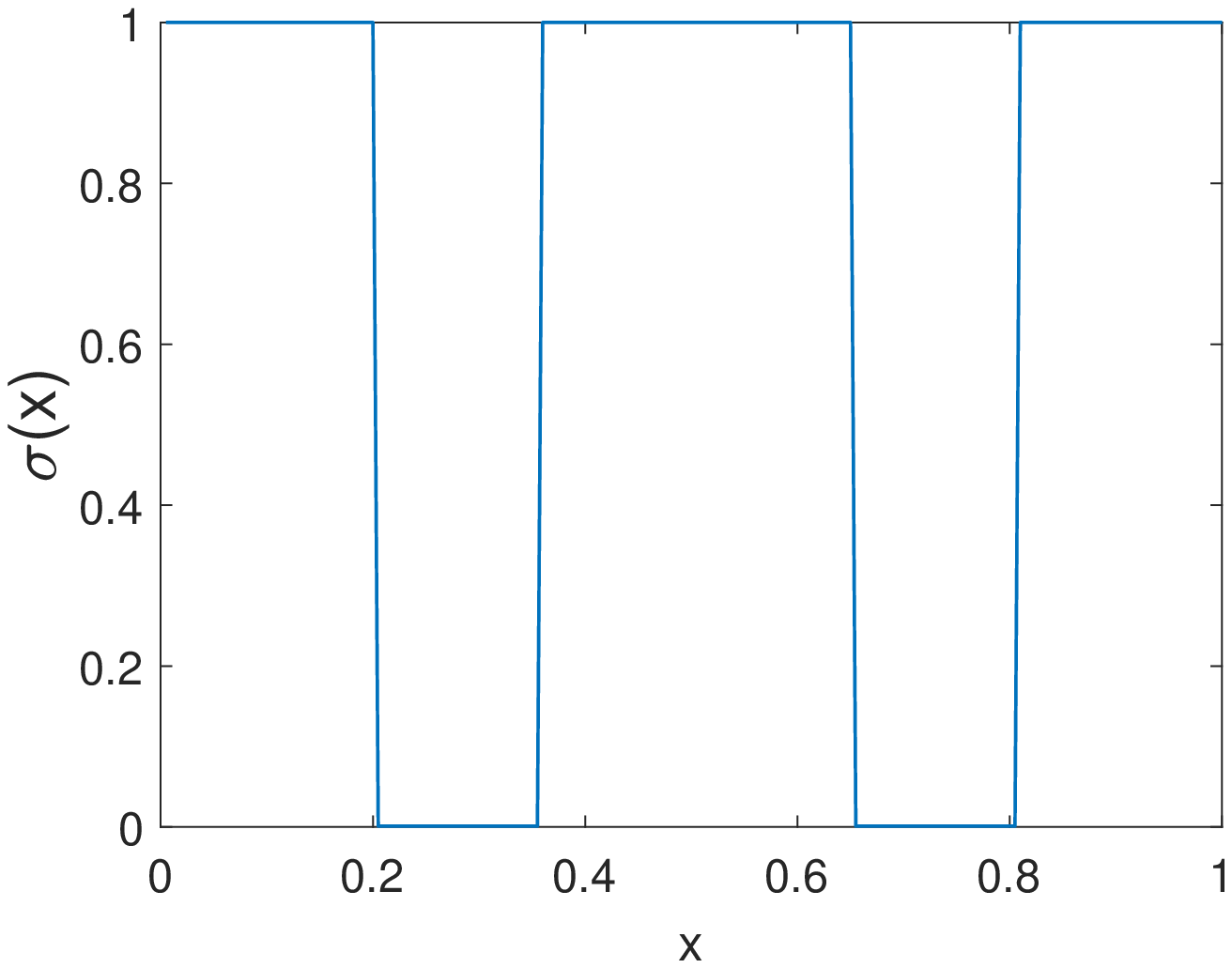}
\includegraphics[width=0.45\textwidth]{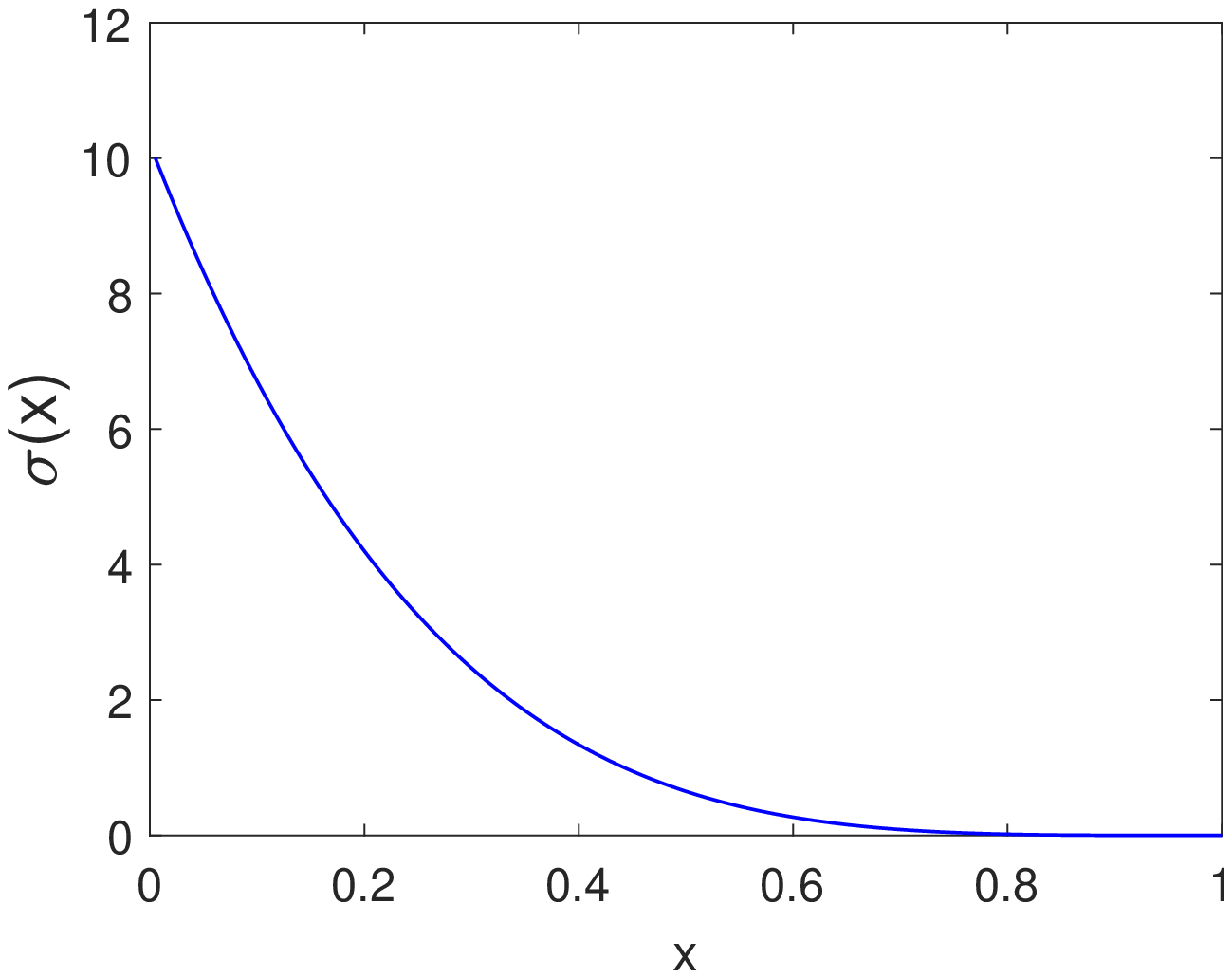}
\caption{Plots of two different collision cross-sections. Left: striped $\sigma$ in \eqref{eqn:sigma}. Right: vanishing $\sigma$ in \eqref{sigma2}. }
\label{fig:sigma}
\end{figure}

\begin{figure}[!h]
\includegraphics[width=0.45\textwidth]{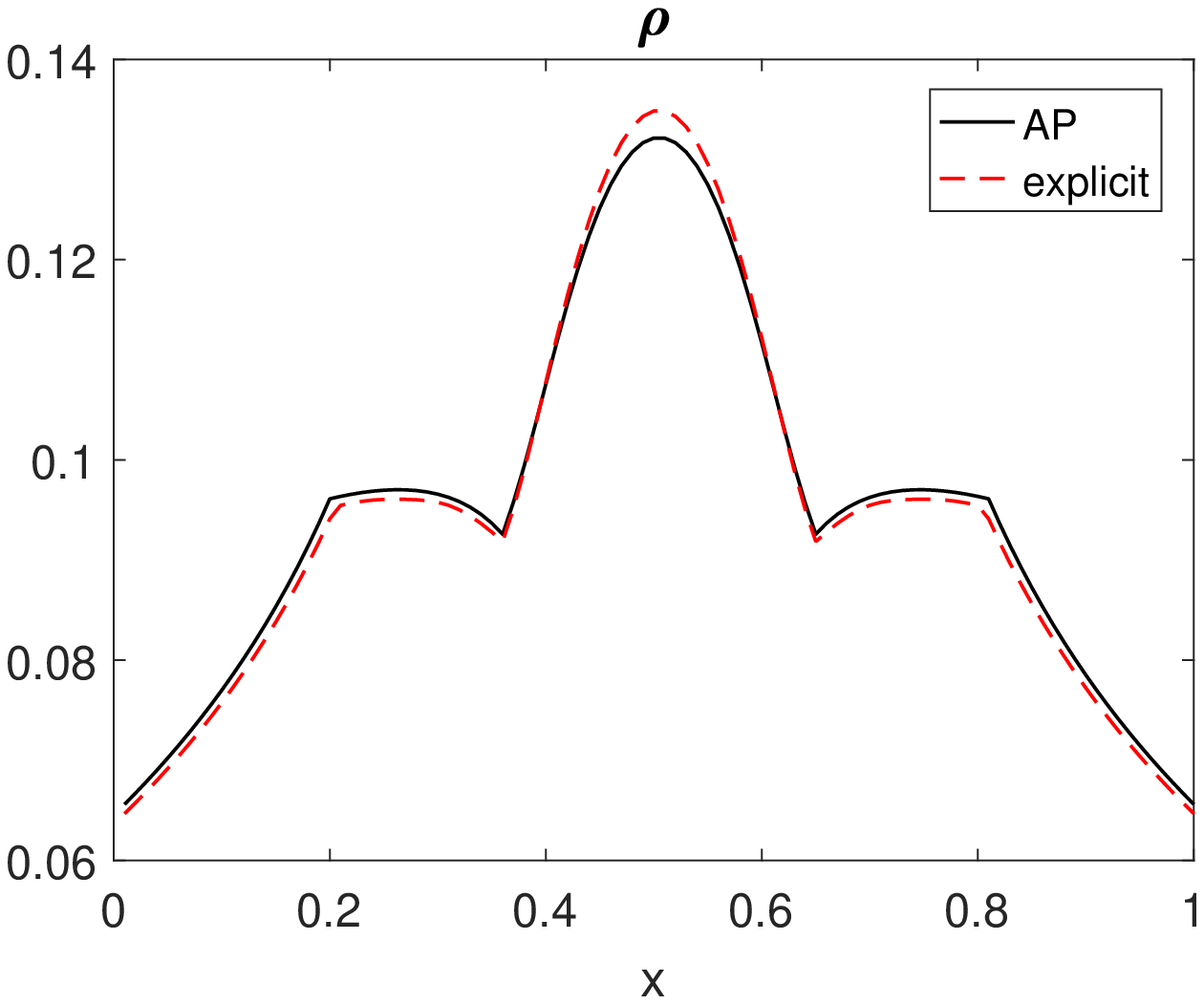}
\includegraphics[width=0.45\textwidth]{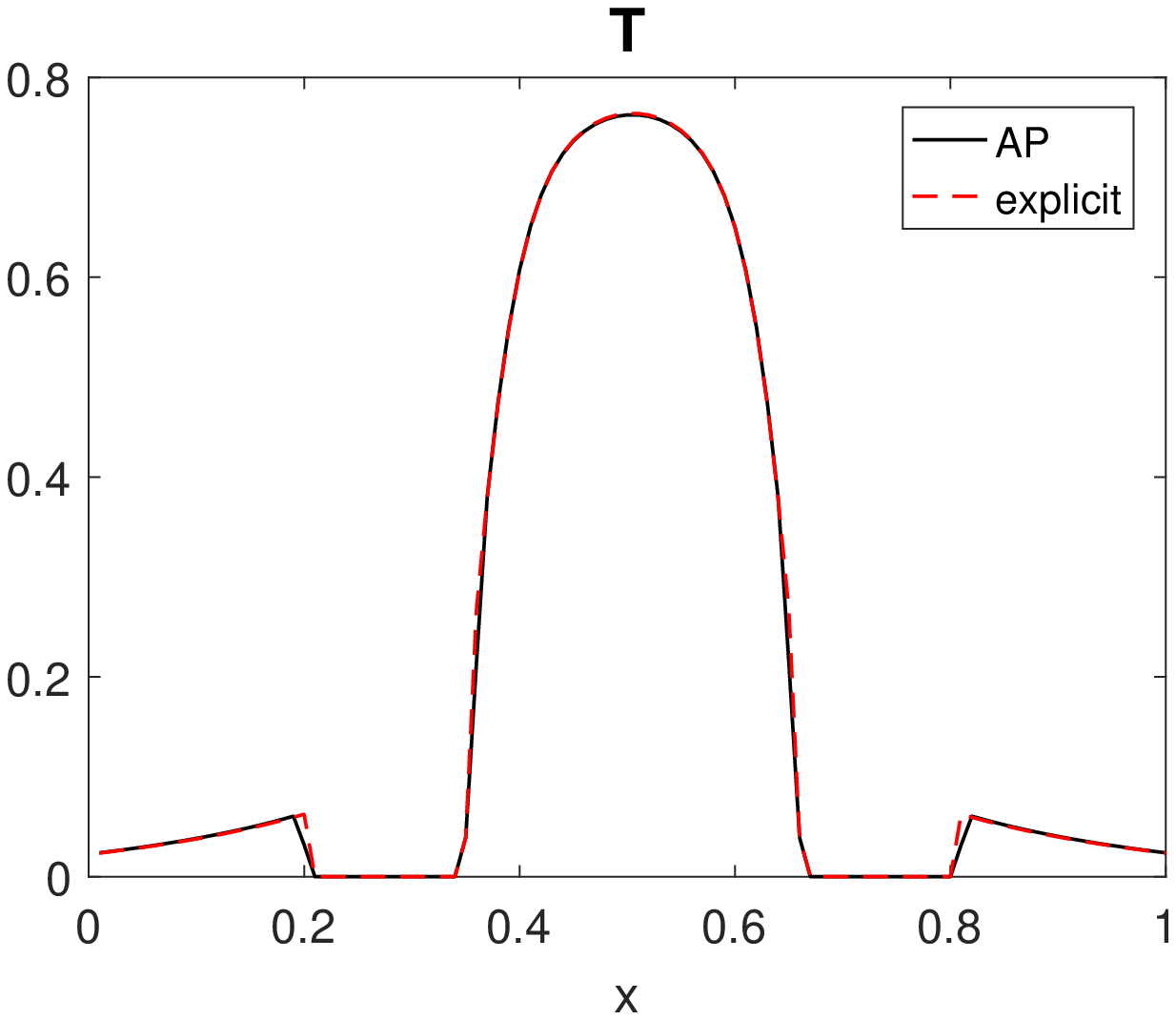}
\caption{Comparison of the density and temperature using our AP scheme and explicit transport scheme at time $t=0.8$. Here $\eps = 1$, $\sigma(x)$ is chosen as in \eqref{eqn:sigma}. Left: density $\rho(x)$. Right: temperature $T(x)$. Here for our AP scheme $\Delta x = 0.01$. For explicit solver, $\Delta x = 0.005$. Time step is chosen as $\Delta t = 0.1 \Delta x$.}
\label{fig:striped_kinetic}
\end{figure}

\begin{figure}[!h]
\includegraphics[width=0.45\textwidth]{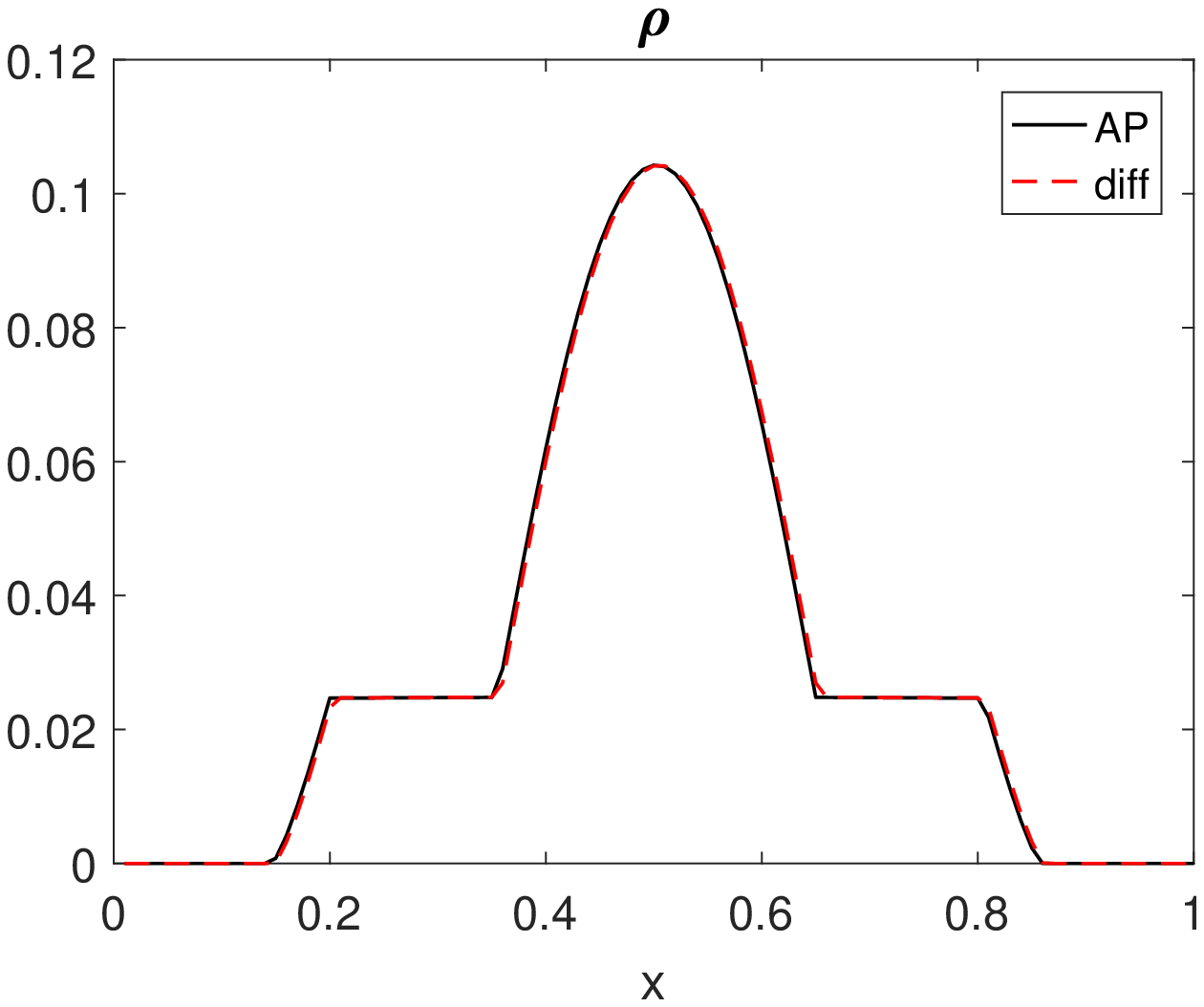}
\includegraphics[width=0.45\textwidth]{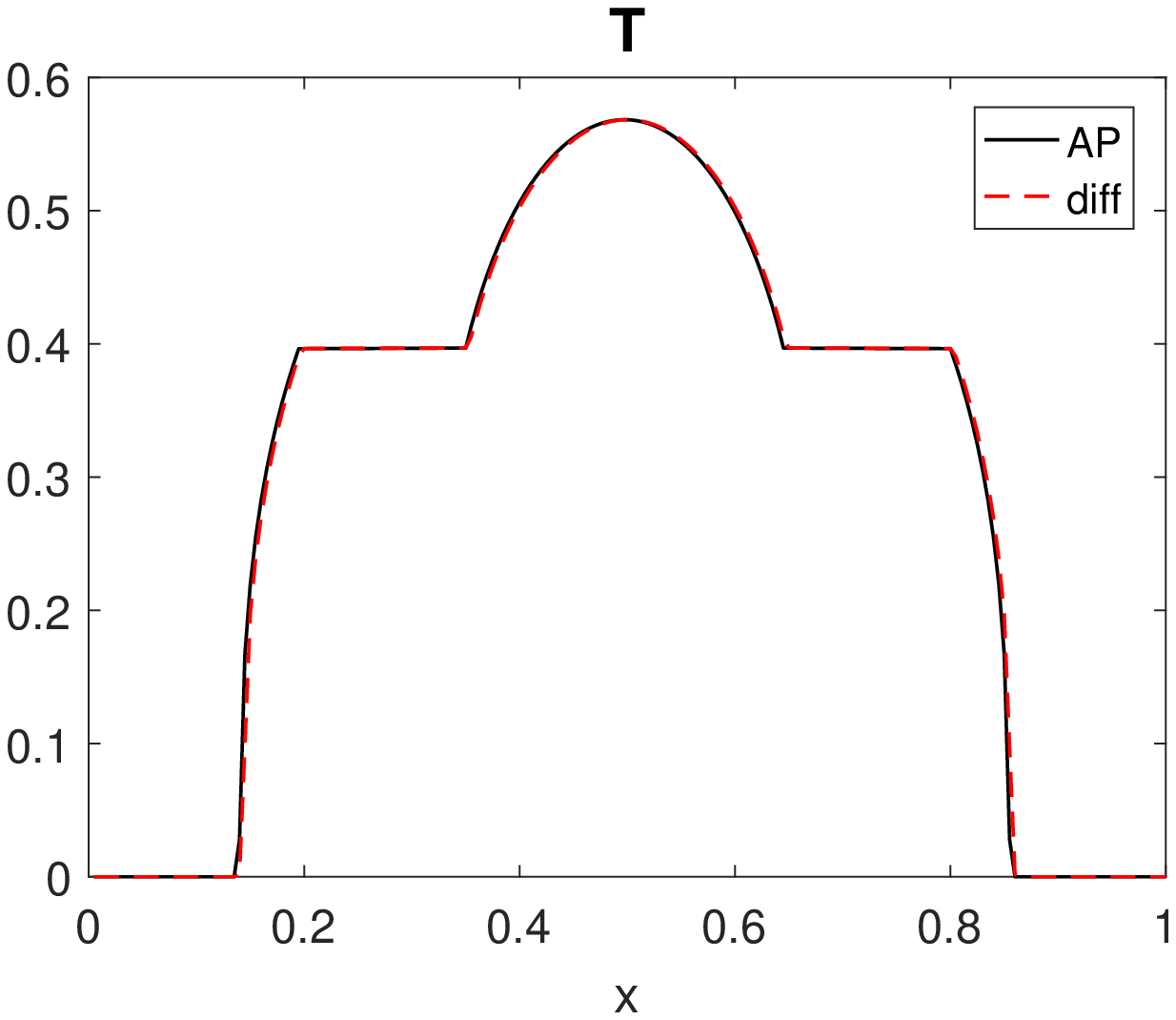}
\caption{Comparison of the density and temperature using our AP scheme with $\eps = 10^{-5}$ and fully implicit solver \eqref{imp_T} at time $t=0.1$. Here $\sigma(x)$ is chosen as in \eqref{eqn:sigma}. Left: density $\rho(x)$. Right: temperature $T(x)$. For both schemes, we use $\Delta x = 0.01$ and $\Delta t = 0.1 \Delta x$.}
\label{fig:striped_diff}
\end{figure}

\begin{figure}[!h]
\includegraphics[width=0.45\textwidth]{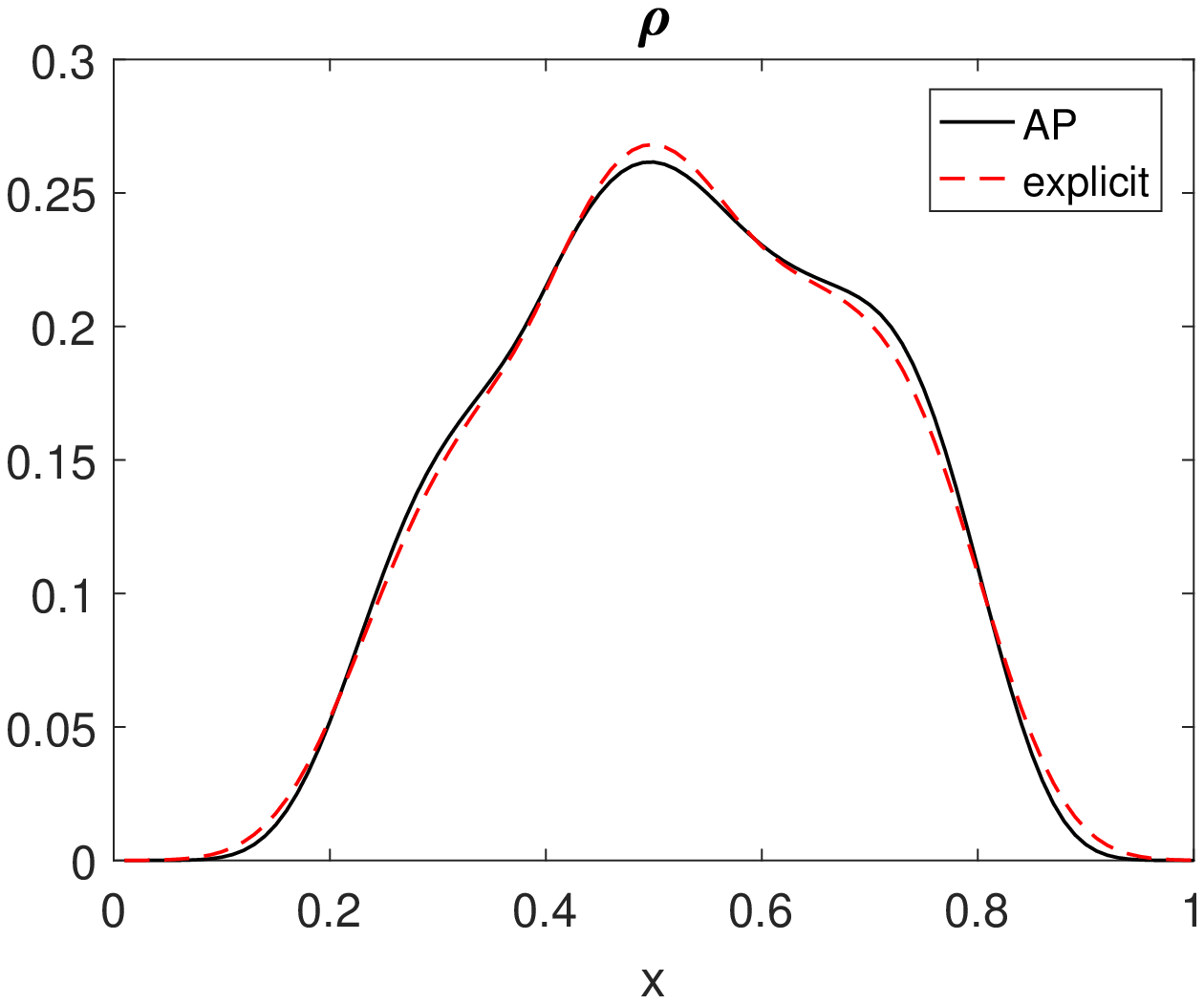}
\includegraphics[width=0.45\textwidth]{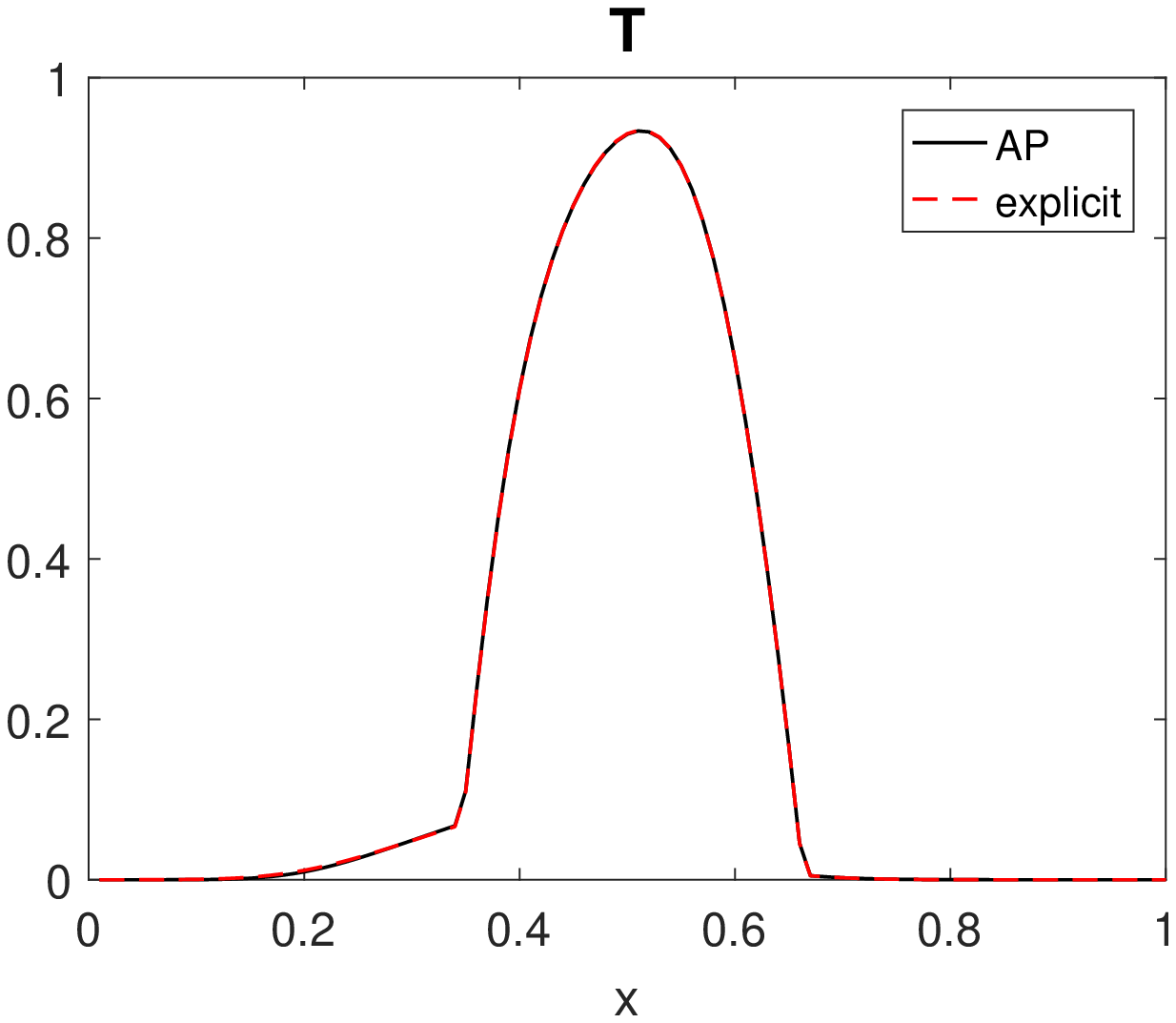}
\caption{Comparison of the density and temperature using our AP scheme and explicit transport scheme at time $t=0.3$. Here $\eps = 1$, $\sigma(x)$ is chosen as in \eqref{sigma2}. Left: density $\rho(x)$. Right: temperature $T(x)$. Here for our AP scheme $\Delta x = 0.01$. For explicit solver, $\Delta x = 0.005$. Time step is chosen as $\Delta t = 0.1 \Delta x$.}
\label{fig:vansig_kinetic}
\end{figure}

\begin{figure}[!h]
\includegraphics[width=0.45\textwidth]{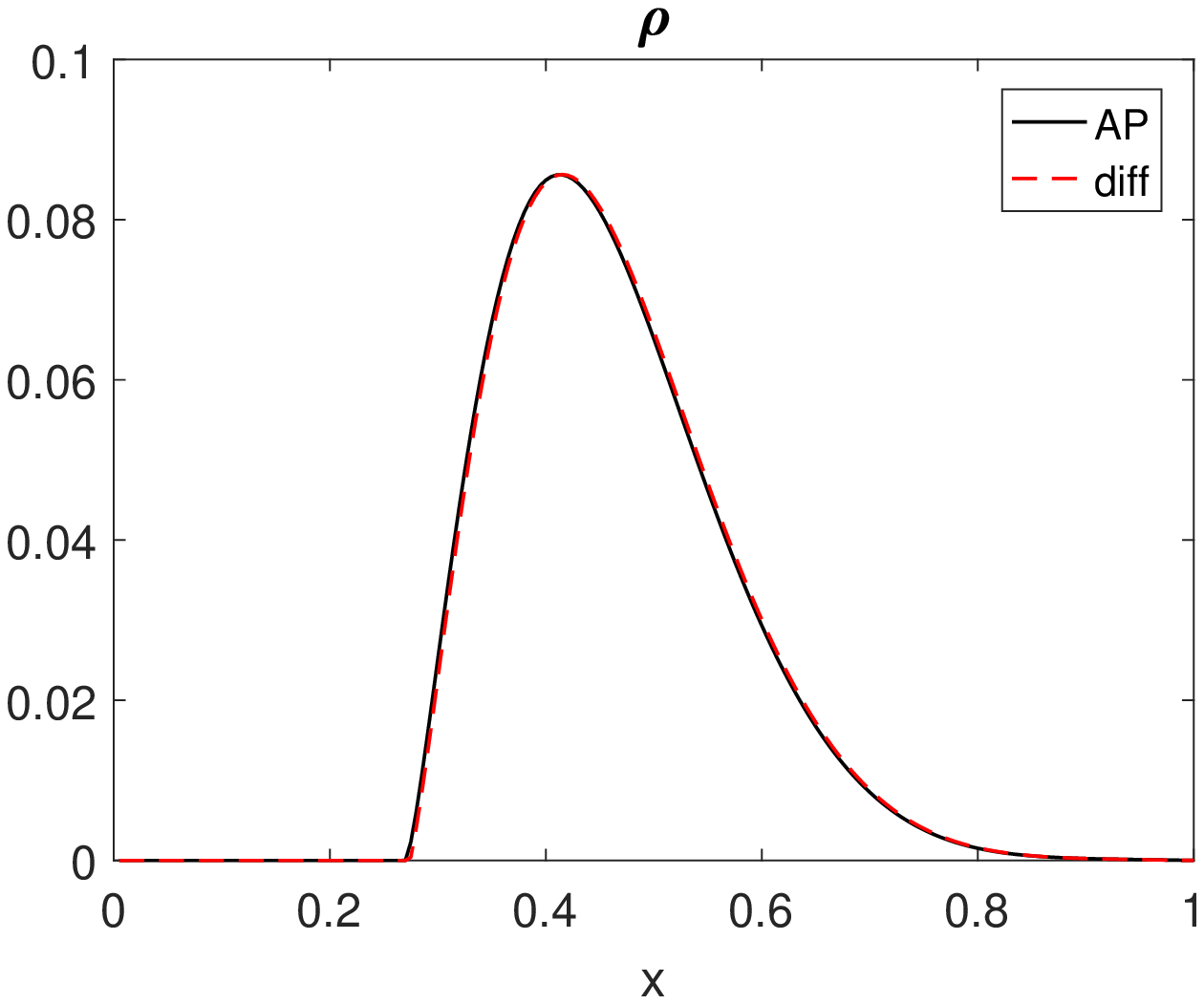}
\includegraphics[width=0.45\textwidth]{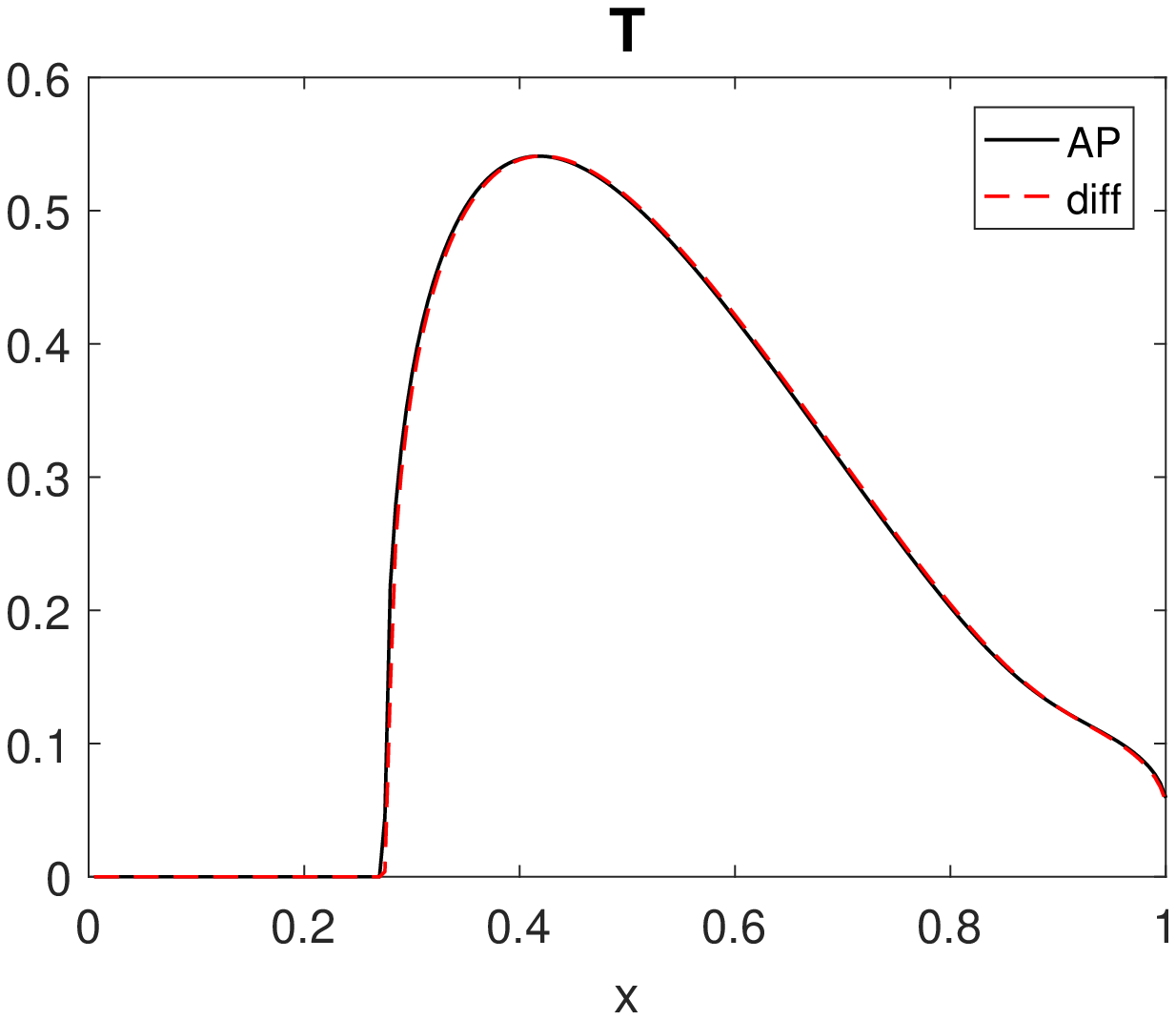}
\caption{Comparison of the density and temperature using our AP scheme with $\eps = 10^{-5}$ and fully implicit solver \eqref{imp_T} at time $t=0.1$. Here $\sigma(x)$ is chosen as in \eqref{sigma2}. Left: density $\rho(x)$. Right: temperature $T(x)$. For both schemes, we use $\Delta x = 0.01$, $\Delta t = 0.1 \Delta x$. }
\label{fig:vansig_diff}
\end{figure}

In this section, we consider two examples with varying $\sigma(x)$. In one case, $\sigma(x)$ is striped and takes the form of  \eqref{eqn:sigma} with $\black{\sigma_0 = 10^{-3}}$. In the other case, $\sigma$ is vanishing in the following sense
\begin{equation} \label{sigma2}
\black{\sigma(x)  = 10(x-1)^4 + 10^{-3}}\,.
\end{equation}
See Fig.~\ref{fig:sigma} for the shape of these two choices of cross-section. Initial data is chosen the same as \eqref{IC000}. In both cases, for $\eps = 1$, we compare the solution to our new scheme with the explicit solver for the GRTE on a finer mesh; for $\eps = 10^{-5}$, we use the solution to the diffusion limit obtained from the fully implicit solver \eqref{imp_T} as a reference. The results are gathered in Figs. \ref{fig:striped_kinetic} -- \ref{fig:vansig_diff}. The good agreement between our solution with the reference solution indicate that our method is capable of dealing with non-smooth cross-section (first case with \eqref{eqn:sigma}), and very week scattering (second case with \eqref{sigma2}) in both optical thin $\eps = 1$ and optically thick $\eps = 10^{-5}$ scenarios.

\subsection{\textcolor{black}{Nonlinear opacity}}
This section consists of two examples with scattering that depends on $T$, i.e., $\sigma(T) = \frac{1}{T^3}$. To compute the reference solution with $\eps = 1$, a direct explicit solver of \eqref{marshak00} would generate erroneous solutions. Therefore, we use an iterative implicit solver for \eqref{pre-marshak}. More specifically, given $I^n$, $U^n$, we run the iteration ($T^{(0)} = T^n$) for $k \geq 0$: 
\begin{equation}  \label{marshak_ref}
\left\{\begin{aligned}
&\frac{1}{c} \frac{I^{(k+1) } - I^n}{\Delta t}+\frac{1}{\varepsilon} v\cdot \nabla_x I^{(k+1)}=
\frac{1}{\varepsilon^2} \frac{\sigma}{{T^{(k)}}^3}\left(acU^{(k+1)}-I^{(k+1)}\right), \quad T^{(k)} = [U^{(k)}]^{1/4} \,;\\
&C_v\frac{ U^{(k+1)} - U^n}{\Delta t}=4\frac{\sigma}{\varepsilon^2}(\rho^{(k+1)}-acU^{(k+1)})
\end{aligned} \right.
\end{equation}
until it converges, and $U^{(k+1)}$, $I^{(k+1)}$ $\rightarrow  U^{n+1}$, $I^{n+1}$. For $\eps \ll 1$, we solve the corresponding diffusion limit \eqref{limit_marshak} via 
\begin{equation}\label{diff_solver}
\frac{(T_i^{n+1})^4 - (T_i^{n})^4}{\Delta t} + \frac{C_v}{a} \frac{T^{n+1}_i - T_i^{n}}{\Delta t} = \frac{4c}{21\Delta x^2} \left[   \frac{ (T_{i+1}^{n+1})^7 - (T_{i}^{n+1})^7 }{\sigma_{i+1/2}}    -  \frac{ (T_{i}^{n+1})^7 - (T_{i-1}^{n+1})^7 }{\sigma _{i-1/2}}   \right]\,,
\end{equation}
where $ \left( \frac{1}{\sigma} \right)_{i+1/2} = \frac{1}{2} \left[  \left( \frac{1}{\sigma} \right)_{i}  +  \left( \frac{1}{\sigma} \right)_{i+1}\right]$.

\begin{figure}[!h]
\includegraphics[width=0.45\textwidth]{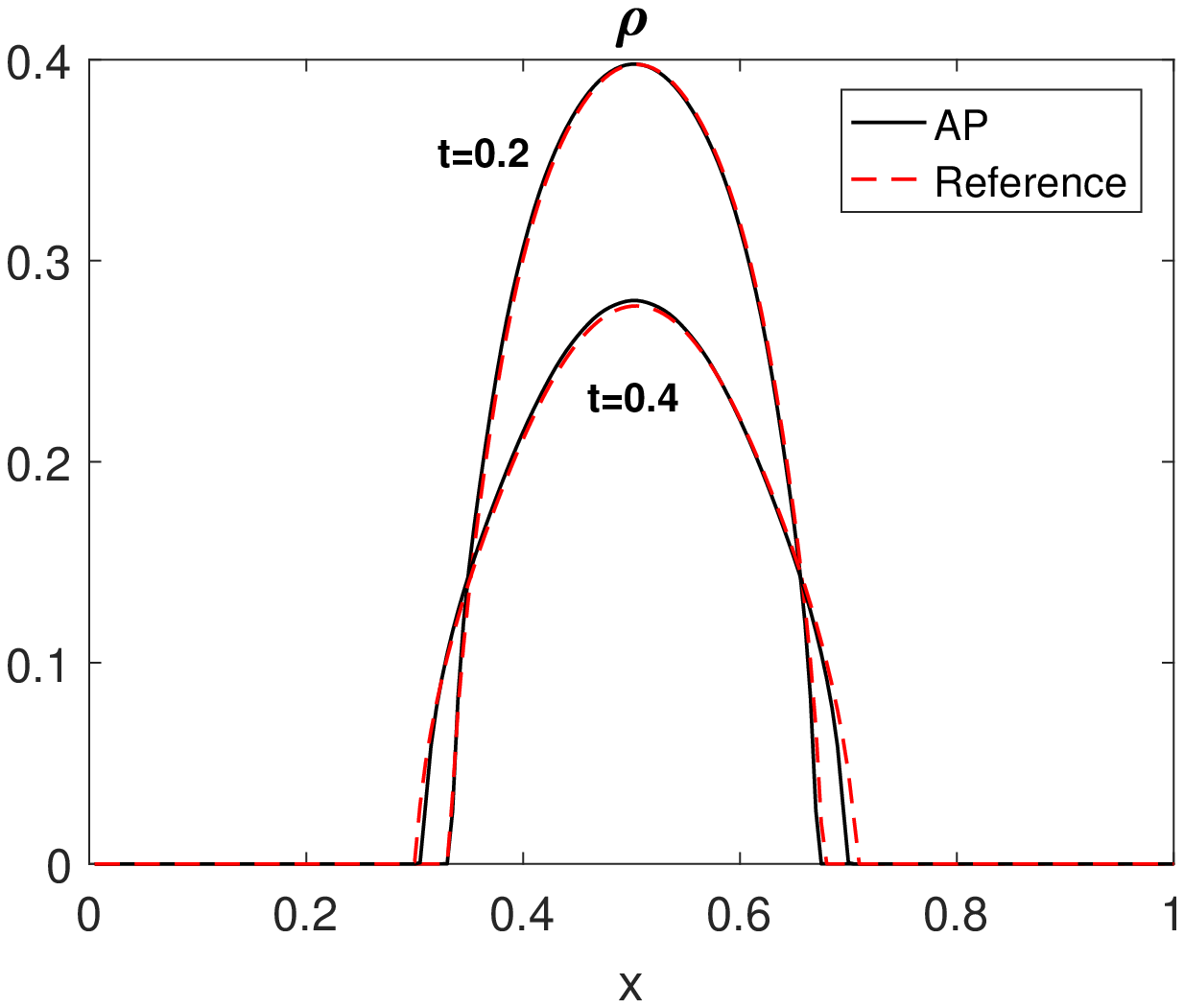}
\includegraphics[width=0.45\textwidth]{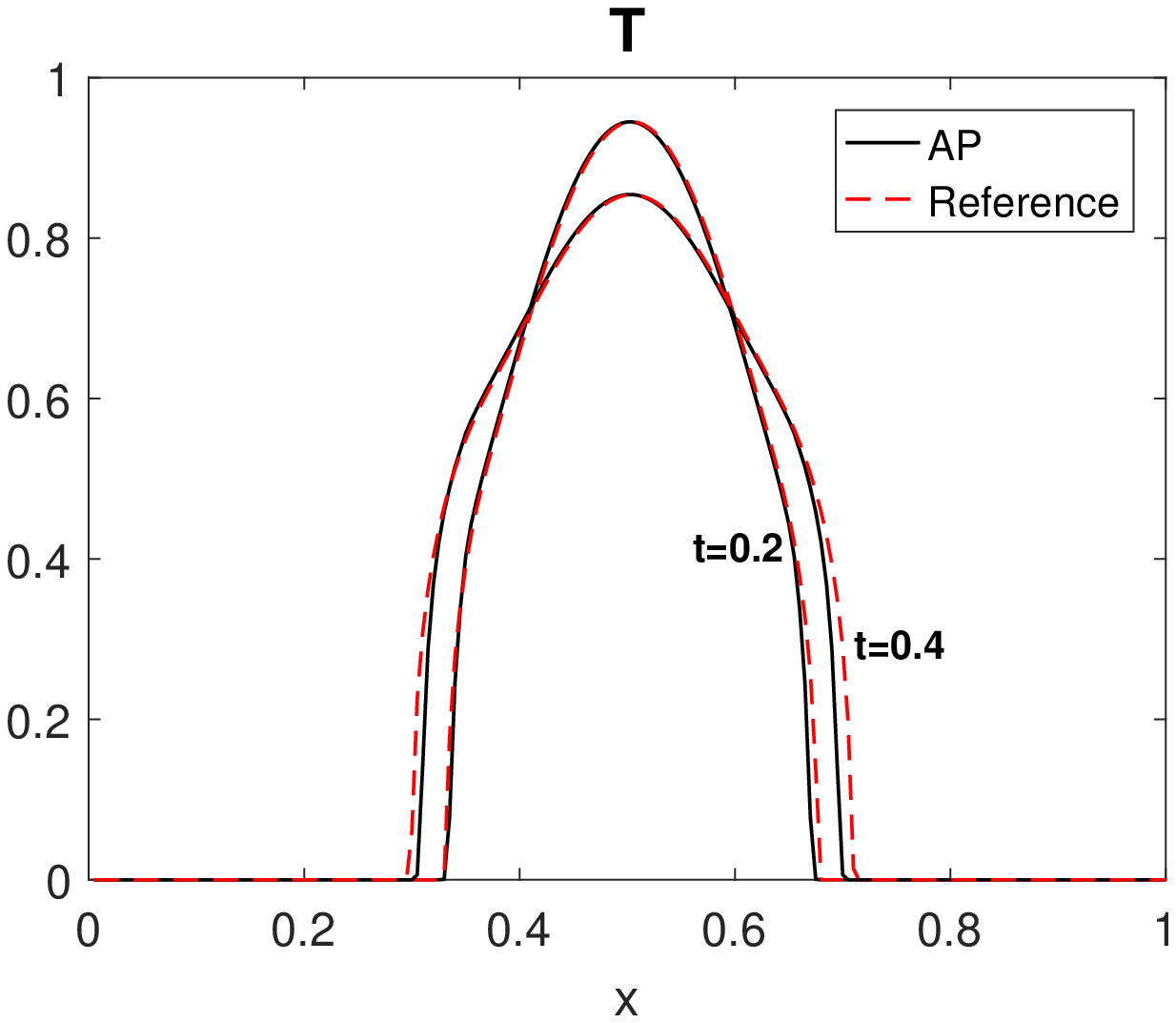}
\caption{Nonlinear opacity case with initial condition \eqref{IC000}, zero boundary condition and $\eps=1$. Comparison of the density and temperature using our AP scheme and reference solver at two different times $t=0.2$ and $t=0.4$. Left: density $\rho(x)$. Right: temperature $T(x)$. Here in our scheme $\Delta x = 0.01$, $\Delta t = 0.1 \Delta x$. For the reference solution, $\Delta x = 0.00125$, $\Delta t  =0.1\Delta x$. }
\label{fig:Marshak1_trans}
\end{figure}

\begin{figure}[!h]
\includegraphics[width=0.45\textwidth]{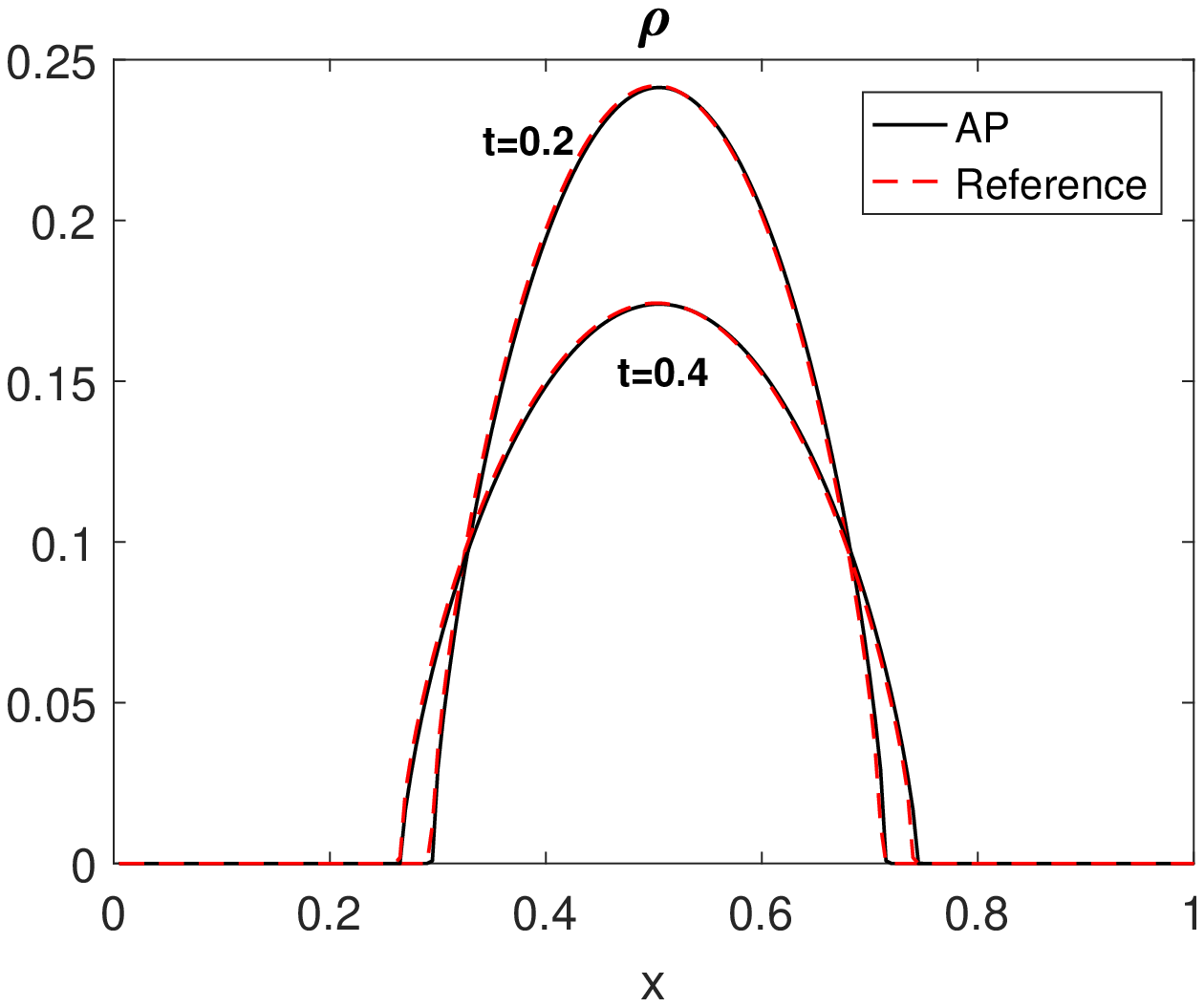}
\includegraphics[width=0.45\textwidth]{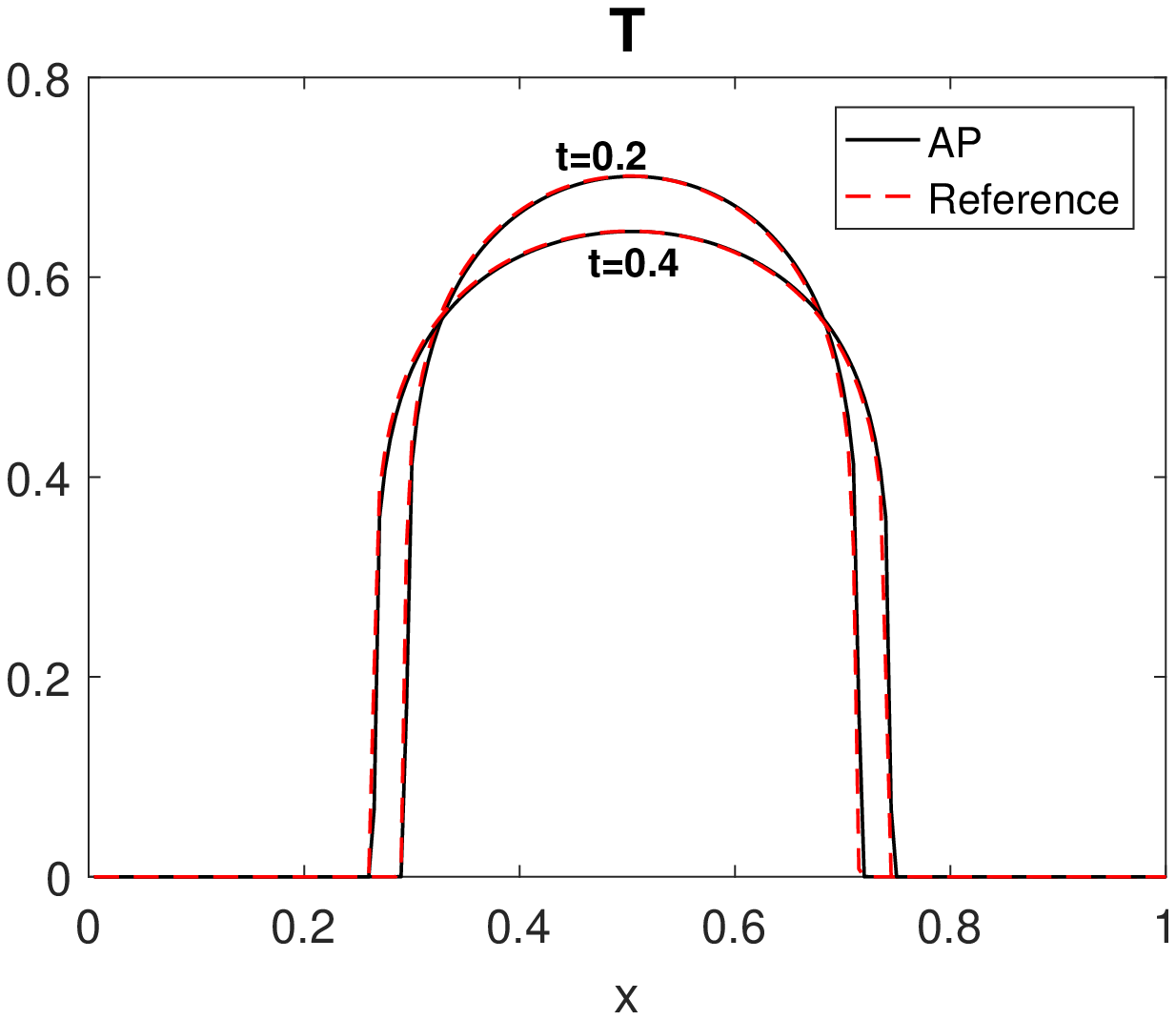}
\caption{Nonlinear opacity case with initial condition \eqref{IC000}, zero boundary condition and $\eps=10^{-5}$. Comparison of the density and temperature using our AP scheme and fully implicit solver \eqref{diff_solver} at two different times $t=0.2$ and $t=0.4$. Left: density $\rho(x)$. Right: temperature $T(x)$. For both schemes we used $\Delta x = 0.01$, $\Delta t = 0.1 \Delta x$. }
\label{fig:Marshak1_diff}
\end{figure}

In the first example, the initial data is taken the same as \eqref{IC000}, and zero incoming boundary condition $I(0,x>0) = I(1,x<0) = 0$ is used. The solutions are plotted in Fig.~\ref{fig:Marshak1_trans} for $\eps = 1$ and Fig.~\ref{fig:Marshak1_diff} for $\eps = 10^{-5}$. 

\begin{figure}[!h]
\centering
\includegraphics[width=0.5\textwidth]{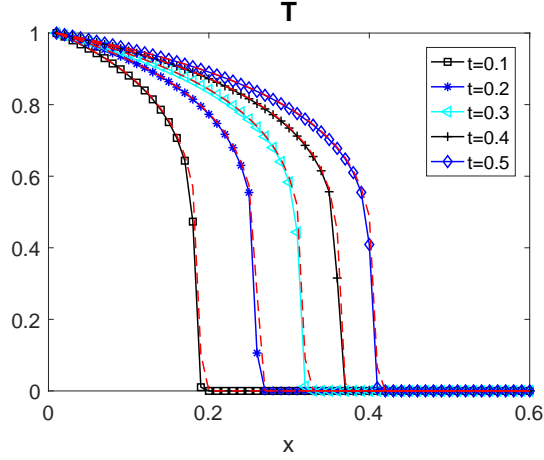}
\caption{Plot of Marshak wave with $\varepsilon = 10^{-5}$ at five different times. The red curves are solutions to the diffusion limit. Here $\Delta x = 0.01$, in our AP scheme $\Delta t = 0.1 \Delta x$, in the diffusion solver \eqref{diff_solver}, $\Delta t = 0.02 \Delta x$. }
\label{fig:Marshak2_1}
\end{figure}

\begin{figure}[!h]
\centering
\includegraphics[width=0.5\textwidth]{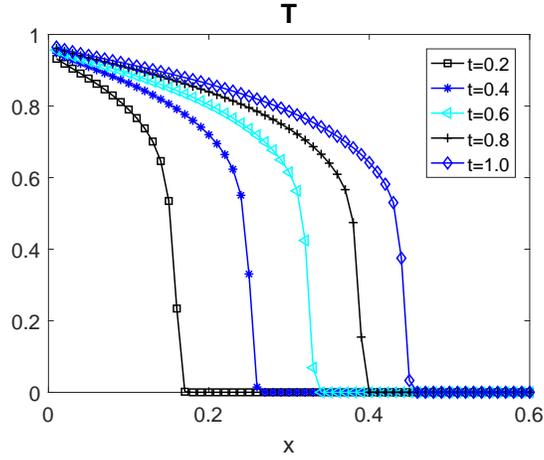}
\caption{Plot of Marshak wave with $\varepsilon = 0.2$ at five different times. Here $\Delta x = 0.01$, $\Delta t = 0.1 \Delta x$. }
\label{fig:Marshak2_2}
\end{figure}

\textcolor{black}{In the second example, the initial data for the transport equation is 
\begin{equation}\label{MarshakIC000}
I(x,v, t=0) = 10^{-16}
\end{equation}
and boundary conditions are
\begin{equation} \label{MarshakBC000}
I(0,v>0) = 1, ~ I(1,v<0) = 0.
\end{equation} 
In the corresponding diffusion limit, the initial temperature takes $T(x,0) = 10^{-4}$, and boundary conditions are $T(x=0) =1 $, $T(x=1) = 0$. In Fig.~\ref{fig:Marshak2_1} we choose $\varepsilon = 10^{-5}$, and compare the solution to our AP scheme with the diffusion solution obtained by \eqref{diff_solver}. We see that our scheme captures the correct propagation front. For moderate $\eps=0.2$, the reference solver \eqref{marshak_ref} fails to produce a reasonable solution, thus we only show the wave propagation obtained from our scheme in Fig.~\ref{fig:Marshak2_2}, which yields a similar wave pattern in the process of time evolution as in the case with much smaller $\eps$. Moreover, we test the stability of this example  with  different $\eps$ and $\Delta x$ in Table \ref{tab:stability3}. Similar as in Table \ref{tab:stability2}, we use the constant $C=\frac{\Delta t}{\Delta x}$ to record the required largest $\Delta t$. We can see that, different from Table \ref{tab:stability2}, $C$ decreases for finer mesh. If the initial data is smooth enough, $C$ becomes uniform for all mesh size, thus the required $C$ depends on the regularity of initial data.}

\begin{table}[!h]    
 \centering
    \fontsize{10}{12}\selectfont    
    \caption{\color{black}Stability test for Marshak wave with initial condition \eqref{MarshakIC000} and boundary condition \eqref{MarshakBC000}.}
\begin{tabular}{l|ccccc}
    \hline
\diagbox[width=5em,trim=l] {$\eps$}{$\Delta x$} & $\f{1}{25}$ & $\f{1}{50}$ & $\f{1}{100}$  & $\f{1}{200}$ & $\f{1}{400}$ \\
\hline
0.2  & 1 & 1 & 1& 1 & 1 \\
1e-03 & 0.5 & 0.3 & 0.1& 0.06 & 0.03 \\
1e-05 & 0.5 & 0.3 & 0.1& 0.06 & 0.03 \\
    \hline
\end{tabular}\vspace{0cm}
   \label{tab:stability3}
    \end{table}
 \subsection{Marshark wave}
\begin{figure}[!h]
\includegraphics[width=0.5\textwidth]{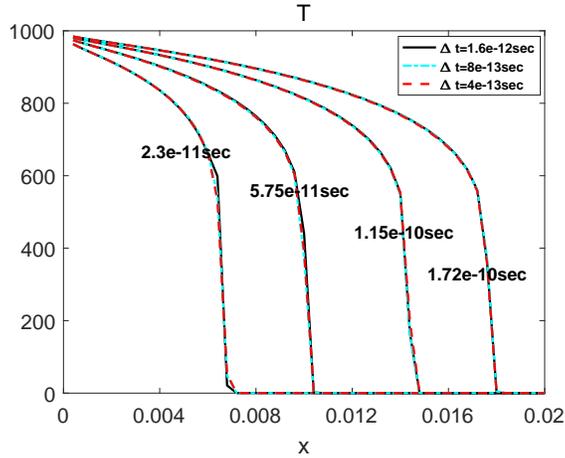}
\caption{\color{black}Marshak wave with units. The initial condition is taken to be \eqref{MarshakIC001}.}
\label{fig:Marshak3}
\end{figure}

\black{
At last, we consider an example with dimensions from \cite{MELD08, LKM13} that reads as follows:
\begin{equation}
\left\{
\begin{aligned}
& \frac{1}{c} \partial_t I+ v\cdot \nabla_x I=
\sigma(T) \left(acT^4-I\right), \quad x  \in \Omega_x \subset \RR^d, ~ v\in \Omega_v \subset \RR^d \,;\\
&  C_v\partial_t T=\sigma(T) (\rho-acT^4), \qquad \rho=\f{1}{|V|}\int_{V}I \dv\,.
\end{aligned}
\right.
\end{equation}
Here the speed of light is $c = 29.98 \text{cm/ns}$, radiation constant $a = 0.01372 \text{GJ}/\text{cm}^3 \text{Ke}V^4$, heat capacity $C_v = 0.3 \text{GJ}/\text{cm}^3\text{KeV}$, and $\sigma(T) = \frac{300}{T^3} \text{cm}^{-1}$. The boundary condition is $T(0,t) = 1$KeV, and initial condition is taken as
\begin{equation} \label{MarshakIC001}
T(x,0) =  \frac{1}{2}[1-\text{tanh}((x-0.0024)*1000)] ~\text{KeV}.
\end{equation}
In Fig.~\ref{fig:Marshak3}, we consider the computational domain $x\in [0,0.02] \text{cm}$, and we plot the solution with different $\Delta t = 1.6\times 10^{-12}$sec, $8\times 10^{-13}$ sec, and $4\times 10^{-13}$sec. It is important to point out that our CFL condition here is comparable to that in \cite{LKM13} and much better than that in \cite{MELD08}. 
}

\begin{color}{black}
\section{Conclusion and discussions}
In this paper, we design a semi-implicit scheme for GRTE that can use a hyperbolic CFL condition independent of $\varepsilon$ to
provide the correct solution behavior. The idea is to introduce a new variable $U$ that takes into account the nonlinearities in the Plank function $T^4$ and then solve a linear system to get a prediction for $T^4$. Then the material temperature $T$ and $I$ are updated by a correction step.

In order to avoid nonlinear iteration, several linearized time discretizations for GRTE have been proposed and studied in the literature \cite{MWS,MELD08}, however, as pointed out in \cite{LKM13}, semi-implicit treatment of the nonlinear Planck function leads to violation of the maximum principle, thus induces oscillations when the time step is not small enough. The time steps have to be chosen as $\Delta t\sim \max\{O(\Delta x^2),O(\varepsilon\Delta x)\}$ in order to achieve the correct solution behavior. Our scheme allows hyperbolic CFL condition $\Delta t\sim\mathcal O(\Delta x)$ with a CFL number independent of $\varepsilon$. 
 The other advantage of our proposed scheme is that it can capture the correct front position for compact supported initial data. Most tests in the literature are for $T$ from very small to very large but not for compact supported data.

 Several problems need further investigations. For example, it remains unclear why the proposed semi-implicit discretization can capture the correct front speed. In fact, this is not always true for arbitrary semi-implicit discretization. In \cite{Lowrie}, the author considered various linearization of nonlinear implicit iterations for the corresponding diffusion equation, and found that some can capture the correct front position while others not. These observations are only numerical and warrants analytical understandings. Another issue is to extend our current scheme to higher dimensions, as we have used even-odd decomposition and staggered grid space discretization, such an extension is not straightforward. These will be our future work.  
\end{color}

\bibliography{nRTE_reference}
\bibliographystyle{siam}

\end{document}